\documentclass{article}
\pdfoutput=1

\usepackage{a4wide}
\usepackage{enumerate}
\usepackage{enumitem}
\usepackage{amssymb,amsmath,amsthm} 
\usepackage{listings}
\usepackage{hyperref}
\usepackage{graphicx}
\usepackage{stmaryrd}
\usepackage{multicol}
\usepackage{tikz}
\usetikzlibrary{shapes,arrows,tikzmark,decorations.pathreplacing,calc}
\usepackage[export]{adjustbox}

\title{Isotopic piecewise affine approximation of\\algebraic or $C^1$ varieties}
\author{Christophe Raffalli}

\newcommand{\hull}[1]{%
  \mathcal H({#1})%
}

\global\let\epsilon\varepsilon

\newcommand{\PDR}{{\cal P}^2(\mathbb{R})}
\newcommand{\PTR}{{\cal P}^3(\mathbb{R})}

\newcommand{\dist}{\mathrm{dist}}

\newtheorem{theo}{Theorem}

\newtheorem{conj}[theo]{Conjecture}

\newtheorem{defi}[theo]{Definition}
\newtheorem{prop}[theo]{Proposition}

\newtheorem{algo}[theo]{Procedure}

\newtheorem{innermanualtheo}{Theorem}
\newenvironment{manualtheo}[1]
  {\renewcommand\theinnermanualtheo{#1}\innermanualtheo}
  {\endinnermanualtheo}

\newtheorem{innermanualconj}{Conjecture}
\newenvironment{manualconj}[1]
  {\renewcommand\theinnermanualconj{#1}\innermanualconj}
  {\endinnermanualconj}

\newcommand{\tr}{{^t\!}}

\newcommand\restr[2]{{#1{\raisebox{-.5ex}{$|$}_{#2}}}\!\!}

\newcommand\reflabel[1]{\stepcounter{equation}\label{#1}(\theequation)}

\makeatletter
\newcommand{\srcsize}{\@setfontsize{\srcsize}{5pt}{5pt}}
\makeatother

\begin{document}

\maketitle

\begin{abstract}
We propose a novel sufficient condition establishing that a piecewise affine
variety has the same topology as a variety of the sphere $\mathbb{S}^n$ defined by
positively homogeneous $C^1$ functions. This covers the case of $C^1$ varieties
in the projective space $\mathbb{P}^n$. We prove that this condition is
sufficient in the case of codimension one and arbitrary dimension. We describe
an implementation working for homogeneous polynomials in arbitrary dimension
and codimension and give experimental evidences that our condition might still
be sufficient in codimension greater than one.
\end{abstract}


\section{Introduction}

\begin{figure}
  \begin{center}
    \includegraphics[width=0.4\linewidth,bb=0 0 1043 1043]{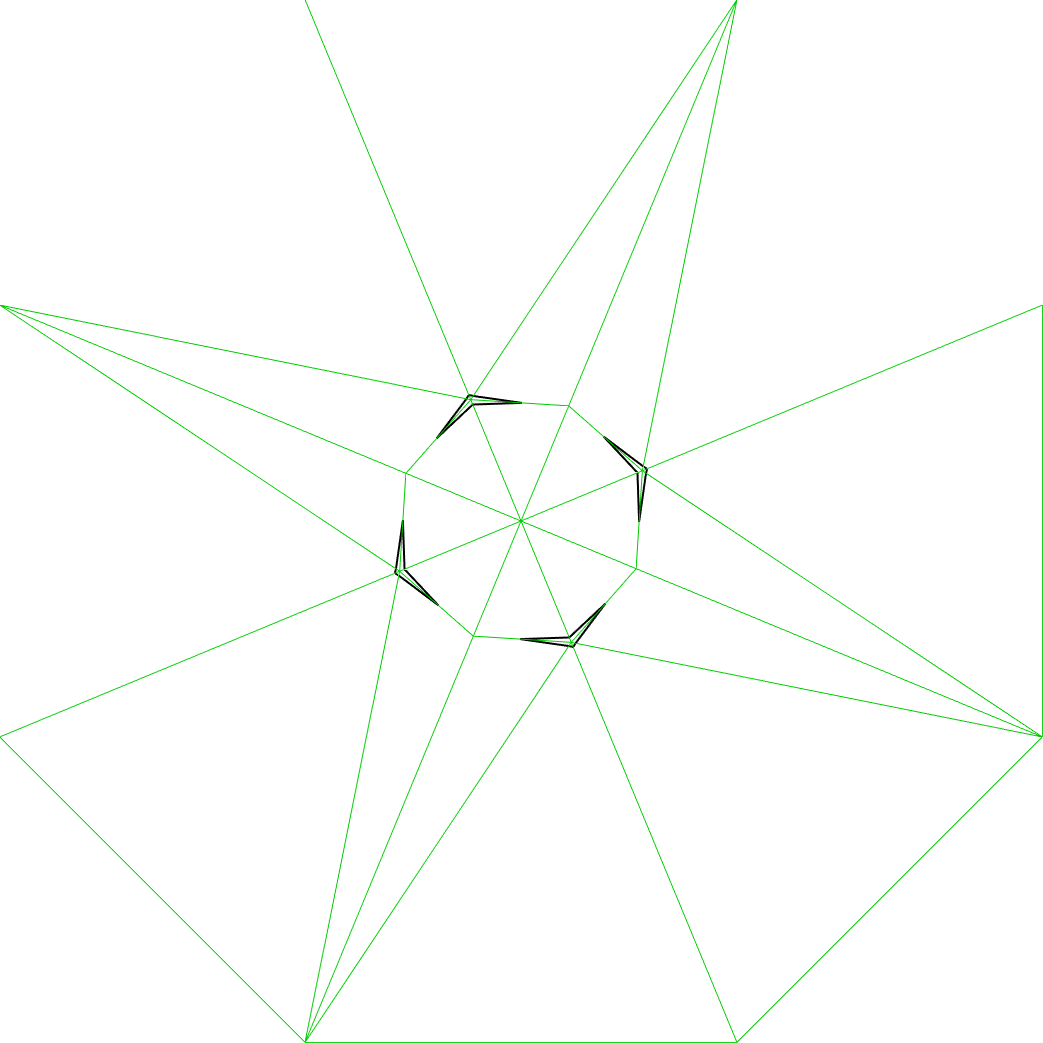}
    \caption{piecewise linear approximation of a quartic curve}\label{figure:quartic_curve}
  \end{center}
\end{figure}

\subsection{Contribution}

Let a variety $V$ be defined by a system of implicit equations: $V = \{x \in
\mathbb{K}, f_1(x) = 0, \dots, f_m(x) = 0\}$ on some compact polyhedron
$\mathbb{K} \subset \mathbb{R}^n$, with $1 \leq m \leq n$.  We assume that the
functions $f_1, \dots, f_m : \mathbb{R}^n \to \mathbb{R}$ are $C^1$.  Let $S =
(S_i)_{i\in I}$ be a decomposition of $\mathbb{K}$ into a family of simplices.
A piecewise affine variety may always be defined from $V$ and $S$ by defining
for each $1 \leq i \leq k$ an approximation $\tilde f_i$ of the function $f_i$
by
\begin{itemize}
 \item $\tilde f_i(x) = f_i(x)$ for any $x$ vertex of $S_i$ for $i \in I$ and
 \item $\tilde f_i$ is affine when restricted to any $S_i$ for $i \in I$.
\end{itemize}
From this, we define $\tilde V = \{x \in \mathbb{K}, \tilde f_1(x) = 0, \dots,
\tilde f_k(x) = 0\}$.  The question is to find a sufficient condition ensuring
that $V$ and $\tilde V$ are isotopic.

Moreover, we search for a criteria that can be computably approximated with
arbitrary precision in the case of multivariate polynomials, to allow for an
implementation.

We propose two theorems in codimension one ($m = 1$) and two
conjectures, one weaker than the other, supported by some experimental evidence,
in the general case. More precisely:
\begin{itemize}
  \item In section \ref{compact}, we give a theorem that answers the question
  when $\mathbb{K}$ is a compact polyhedron in $\mathbb{R}^n$, in codimension one
  ($m = 1$) and when $f_1$ is of $C^1$ class.
  \item In section \ref{projective}, we show that the same condition is
    correct if $\mathbb{K} = \mathbb{S}^n$ the unit sphere of
    $\mathbb{R}^{n+1}$, in codimension one and when $f_1$ is of $C^1$
    class and positively homogeneous of degree $d$ (i.e. $f(\lambda x) = \lambda^d
    f(x)$ for all $\lambda \in \mathbb{R}_+$ and $x \in \mathbb{R}^{n+1}$).
    This case could be considered as codimension $2$, but the homogeneity allows
    to ignore completely the equation of the sphere.
  \item In section \ref{codimension}, we generalise the previous statement to
    arbitrary codimension and conjecture that it still holds.  This conjecture
    is supported by an implementation which was tested on many examples and
    never gave wrong topology. We actually give two conjectures, because we
    lack some results on convex set of full rank matrices. We give more details
    in subsection \ref{another} and section \ref{codimension}.
\end{itemize}

Let give now the statement of our first theorem in section \ref{compact}:
\begin{manualtheo}{\ref{thcompact} page \pageref{thcompact}}
  Let $\mathbb{K} \subset \mathbb{R}^n$ be a compact polyhedron.
  Let $(S_i)_{i\in I}$ be a simplicial decomposition of $\mathbb{K}$.
  Let $f : \mathbb{R}^n \to \mathbb{K}$ be a $C^1$ function in $n$ variables.
  Let $V = \{ x \in \mathbb{K}, f(x) = 0\}$ be the zero locus of $f$
  restricted to $\mathbb{K}$.
  Assume that $V \cap \partial\mathbb{K} = \emptyset$.

  We define $\tilde{p} : \mathbb{K} \to \mathbb{R}$ the piecewise affine
  function such that for all $i \in I$, $\restr{\tilde f}{S_i}$ is affine
  and for any $v$ vertex of $S_i$, we have
  $f(v) = \restr{\tilde f}{S_i}(v)$.
  We define the following:
  \begin{itemize}
  \item $\tilde V = \{ x \in \mathbb{K}, \tilde f(x) = 0\}$ the zero locus of
  $\tilde f$.
  \item $\mathbb{K}(f) = \{ x \in \mathbb{K}, f(x)\tilde{f}(x) \leq 0\}$.
  \item $\tilde \nabla f(x) = \{ \nabla \restr{\tilde f}{S_i}(x), x \in S_i
    \} \subset \mathbb{R}^n$
  \item $G(f,x) = \{ \nabla f(x) \} \cup \tilde \nabla f(x) \subset \mathbb{R}^n$
  \end{itemize}
  If the condition (\ref{ieq1}) below holds, then $V$ and $\tilde V$ are isotopic:
  \begin{eqnarray}
    \forall x \in \mathbb{K}(f), 0 \notin \hull{G(f,x)} \text{ the convex
      hull of $G(f,x)$}\label{ieq1}
  \end{eqnarray}
\end{manualtheo}

Let us give some ideas about this theorem: the isotopy is naturally defined by
$f_t(x) = t f(x) + (1 - t) \tilde f(x)$ for $x \in \mathbb{K}$ and $t \in
[0,1]$. The function $x \mapsto f_t(x)$ is not differentiable, but it is
differentiable in any direction, and its gradient at $x$ in direction $D$ is
always given by a scalar product $V.D$ where $V$ is in the convex hull of
exactly the gradients we are considering in the set $G(f,x)$. Thus our
condition ensures that $V \neq 0$, which we find is very natural
\emph{smoothness} condition. This condition only needs to hold in region where
$f_t$ may be null, i.e. when $f(x)$ and $\tilde f(x)$ have opposite sign, this
justifies the definition of $\mathbb{K}(f)$.

The theorem of section \ref{projective} is almost the same, we ask
for the function to be positively homegenous and we decompose
$\mathbb{R}^{n+1}$ in \emph{simplicial cones}, which is defined in section
\ref{convention}. Appart from this, the statement is unchanged.

In section \ref{codimension} we propose two conditions (conjecture \ref{conj}
and \ref{wconj}) that could apply in arbitrary codimension (i.e. with more than
one polynomials). Unfortunately we are not able to prove those conjectures.

The first one (conjecture \ref{conj}), the most natural, generalises the
condition (\ref{ieq1})
$$
  \forall x \in \mathbb{K}(p), 0 \notin \hull{G(p,x)}
$$
into
\begin{eqnarray}
  \forall x \in \mathbb{K}(p), \forall A \in \hull{G(p,x)}, A \text{ is of
    maximal rank} \label{ieq2}
\end{eqnarray}

This is natural as with codimension greater than one, the gradients become
matrices and maximal rank expresses transversality, hence smoothness.

Remarks: we do not need extra hypotheses, like smoothness (or non complete
intersection with codimension greater than one). However, if the
variety is not sm
ooth, our criteria will never be satisfied. We should also
say that our condition is frame independant. Indeed, a change of
coordinates will multiply all the gradients by the same invertible matrix and
the convex hull is transformed accordingly.

\subsection{A global criteria}

A standard way to compute piecewise affine approximation of varieties defined by
implicit equations are decomposition method that proceed by incrementally
subdivising the ambient space in simplices or hypercubes, until some criteria is met.

Our criteria is such a stopping condition, but it is global. To our knowledge,
all existing criteria (like in \cite{LMP07}) will ensure the isotopy of the
orginal variety and its approximation when restricted to each simplex of
hypercube. This is not the case of our criteria.

Let us explain more precisely what we mean by \emph{global}. From the
definitions in our theorem or conjecture, if follows that if $X$ is the set of
vertices of the simplicial decomposition, we have $X \times V$ isotopic to $X
\times \tilde V$. This means that the isotopy can not traverse the vertices of
the decomposition. However it may traverse faces of simplices of dimension one
or more, allowing to use less simplices.

\begin{figure}[ht]
  \begin{center}
    \includegraphics[width=0.6\linewidth,bb=0 0 582 411]{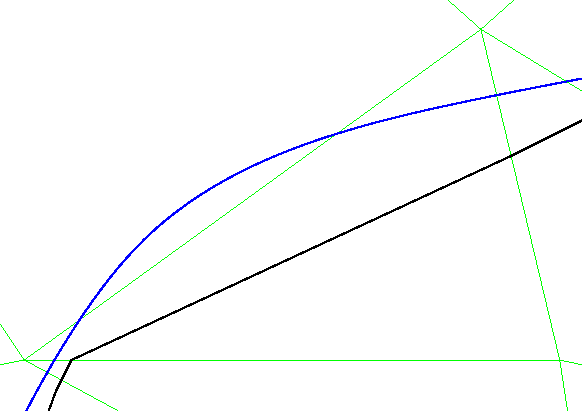}
    \caption{Our criteria is global}\label{figure:global}
  \end{center}
\end{figure}

This is illustrated by figure \ref{figure:global}. This figure represents a
piece of a $C^1$ variety in blue and
its approximation in black. In the triangle which is fully displayed, the
approximation has only one component while the original variety, has two
components. Still our criteria accepts this decomposition.
This can save quite a lot of triangles.

\subsection{Testing the criteria}

Implementing a test for the criteria is not possible in general. But, in the
case of polynomials, we can use Bernstein basis: for all $x$ in a simplex $S$,
$\nabla p(x)$ always lies in the convex hull of the coefficients of the
polynomial $\nabla p$, in this basis, after a change of variable to send the
unit simplex in $S$. This gives easily a sufficient condition to satisfy the
test in each face of the simplicial decomposition. This is detailed in section \ref{implementation}.

Moreover, we can approximate the real criteria with arbitrary precision by
subdividing each face to test. We only do this subdivision for the
test. Refining the decomposition requires to
subdivise the neighbour simplices and we do not want to do that if we can
avoid it.

We implemented a heuristic that searches for a simplicial decomposition
satisfying our criteria. This is relatively quick because it uses
floating point arithmetic. But when an apparently correct decomposition is
found, we retest the criteria with exact rational arithmetic, ensuring the
correctness.

Moreover, the search for the decomposition produces certificates that the
relevant convex hulls do not contain $0$. In codimension one, such a
certificate is a vector which has a positive scalar product with all the
generators of each convex hull. This way the only computation we have to do in
exact arithmetic are change of coordinates, scalar products and comparison. We
do not perform nested computation in loops and this limits the growth of the
size of numerators and denominators. In our experiments, the computing time of
the final exact test is faster than the search for a simplicial decomposition.

\subsection{Examples}

\paragraph*{Remark}:
Because polynomials have no well defined value in the projective space, we
will work within $\mathbb{R}^{n+1}$ and its unit sphere $\mathbb{S}^n$. This is
equivalent and much easier. Still all the examples will be depicted in the
projective space as it avoids to draw every point  of the variety twice.

\begin{figure}
  \begin{center}
    \includegraphics[width=0.4\linewidth,bb=0 0 978 933]{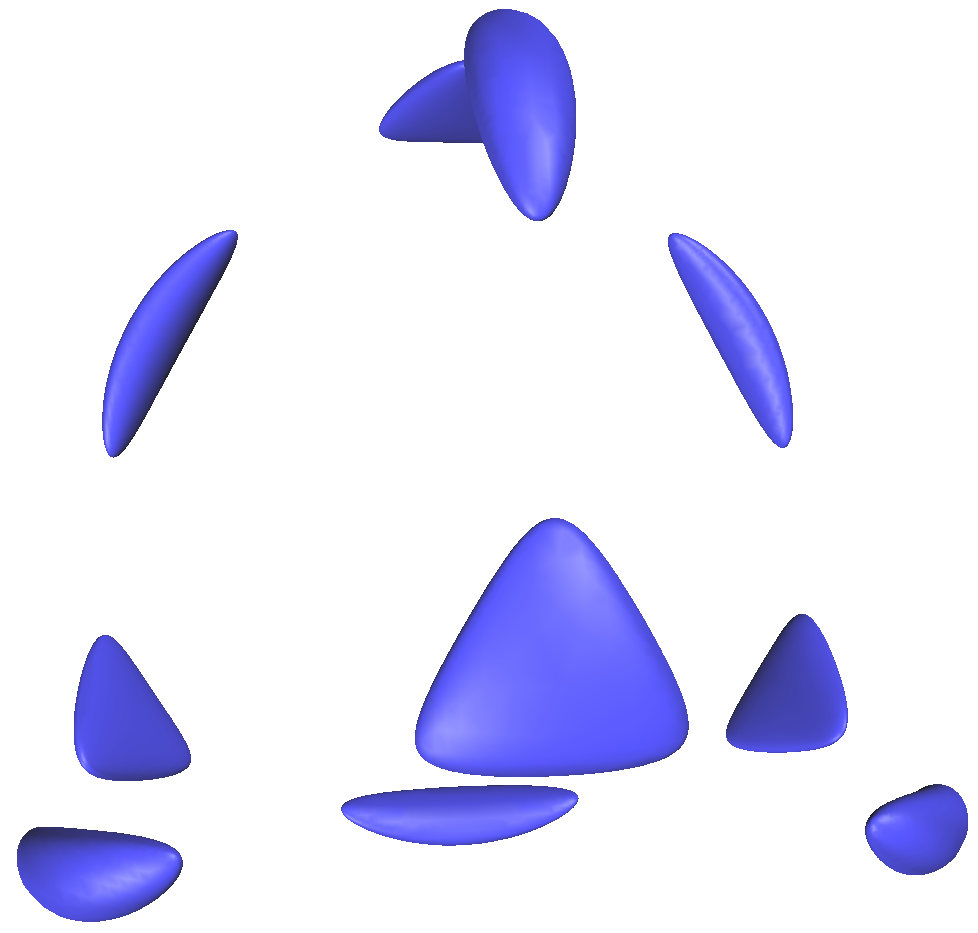}
    \hspace{1cm}
    \includegraphics[width=0.4\linewidth,bb=0 0 1058 940]{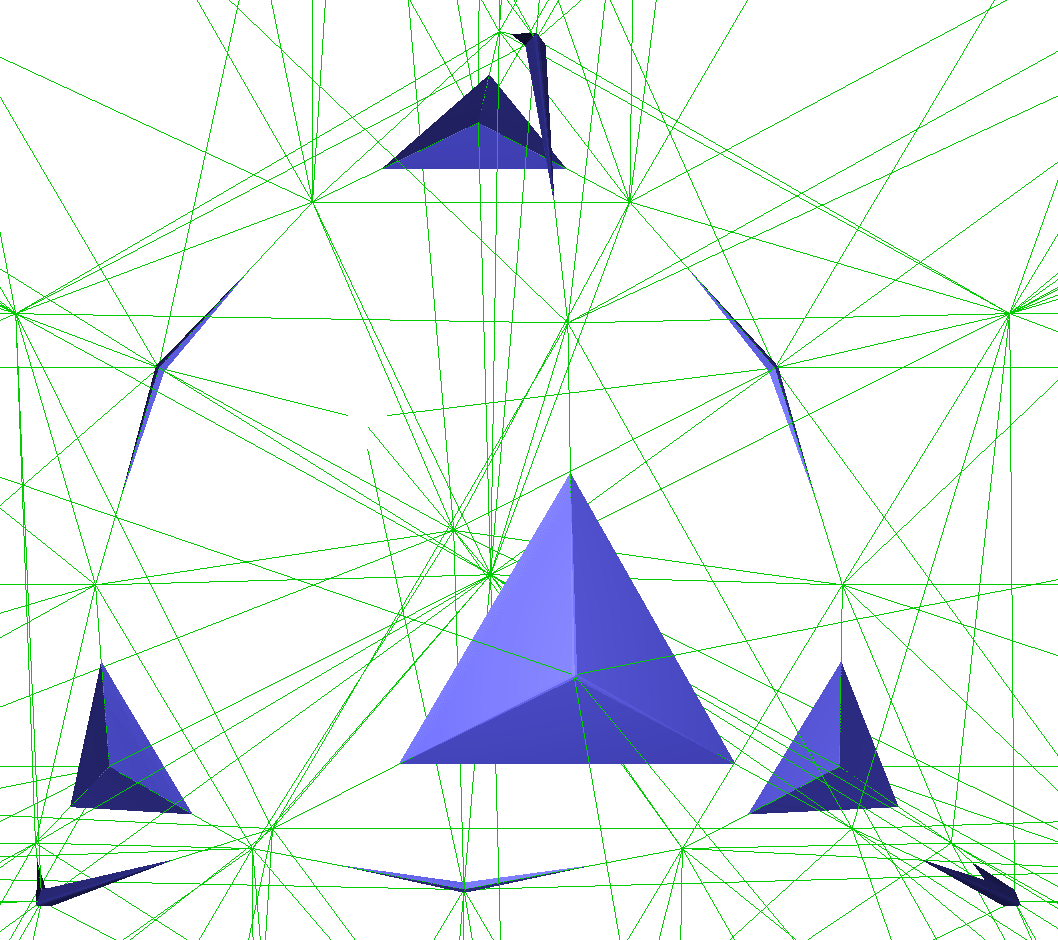}
    \caption{a quartic surface and its piecewise affine approximation}\label{figure:quartic}
  \end{center}
\end{figure}

The green line segments in figures \ref{figure:quartic_curve} and \ref{figure:quartic} are the edges of a simplicial
decomposition of $\PDR$ and $\PTR$ respectively (in the latter case, it is
unfortunately not
easy to guess the simplices from their edges).
The figure \ref{figure:quartic_curve} gives the piecewise affine
approximation of a plane curve of degree 4. It uses a decomposition of the
projective plane with 13 vertices and 24 triangles and requires 58ms to compute. The figure \ref{figure:quartic}
shows an algebraic surface of the same degree together with its
approximation. The decomposition uses 32 vertices and 152 tetrahedron and
requires 3.3s to compute.

As those varieties are enclosed in a compact polyhedron, their approximations
can be proved isotopic to the original varieties by the previous theorem
\ref{thcompact}. They can also be proved correct by the theorem
\ref{thprojective} of section
\ref{projective}.

\subsection{Another conjecture}\label{another}

A key ingredient for both the proof and the implementation is the geometric
form of Hahn-Banach theorem: in codimension 1, if $S$ is a finite set of
vectors, from $0 \notin \hull{S}$, we get a vector $N$ such that $N.V > 0$ for
all $V \in S$. Unfortunately, we could not find nor prove a similar result for
convex sets of full rank matrices. This suggests the following very
interesting conjecture:

\begin{manualconj}{\ref{matconj}}
  Let $1 < m \leq n$ be two natural numbers, let $S \subset
  \mathrm{Mat}_{m,n}(\mathbb R)$ be a
  convex set of matrices of rank $m$. There exists a matrix $M \in
  \mathrm{Mat}_{m,n}(\mathbb R)$ such that $M \tr A + A \tr M$ is symmetric definite and
  positive for all $A \in S$.
\end{manualconj}

This Hahn-Banch conjecture also allows for a notion of \emph{certificate}
allowing to search for a decomposition using efficient floating point
computations, and rechecking the final result using exact rational arithmetic.
We can check quickly that $M \tr A + A \tr M$ is symmetric definite and
positive using Choleski decompostion for each matrices $A$ in $S$ (if $S$ is
finite). This does not require to compute the spectrum of the matrix.

Unfortunately, we do npt know how to implement a test to decide if all
matrices in a convex set of matrices are fullrank (used in \ref{ieq2}), that produces
a certificate. A constructive proof of conjecture \ref{matconj} would likely
provide such a test.

To be able to propose an implementation working in codimension greater than
one, we use a simpler sufficient condition, using only scalar products, that
implies (\ref{ieq2}). This means we have stronger evidence for another
conjecture \ref{wconj} using this stronger condition than for conjecture
\ref{conj} using (\ref{ieq2}).

\subsection{Search for a simplicial decomposition}

Our implementation is described in section \ref{implementation} and is
available on github:
\begin{quote}
  \verb!https://github.com/craff/hypersurfaces!.
\end{quote}
It implements a semi-algorithm building a simplicial decomposition satisfying
our criteria. This means we can solve non-degenerate systems of homogeneous
real polynomial equations. Non-degenerate means that the jacobian matrix of the
system is full rank at point that are in the solution.  By \emph{solving} we
mean finding a piecewise affine approximation of the solution that is isotopic
to the real solution. It allows to compute topological invariants of the
solution and in particular the number of connected components and the Betti
numbers of each component (see for instance the appendix B.3 of
\cite{mangolte2020real} for definitions). It is a semi-algorithm, which means
that it may loop when the system of polynomials is degenerated. Our
semi-algorithm terminates, in principle, if we know that the system of polynomials is
non-degenerated.

In codimension one, our implementation is exact and provides a proved result
about the hypersurface that reposes only on the correction of the final test
using the certificate. This test is a rather short piece of code and would
be easy to rewrite. For codimension greater than one, the validity of the
answer of our algorithm depends upon conjecture \ref{wconj}.

The search for an adequate simplicial decomposition is in an early
stage. Still some examples are already quite interesting. For instance, we are
able to compute Betti numbers of random cubic given by a random homogeneous
polynomials up to dimension 6 in less than one hour. Moreover, to our
knowledge, there is no available exact implementation working in arbitrary
dimension and codimension.

\subsection{Related work}

\subsubsection{Viro's method}

Our main result may be seen as some inverse of Viro's method
\cite{Vir83,Ris93}, more precisely the combinatorial patchworking version for
complete intersection \cite{Stu94}. Indeed, Viro's method allows for the
construction of polynomials with a zero locus having the same topology than a
given piecewise linear variety. We do the opposite.

An important difference is that we use arbitrary simplicial decomposition
rather than the Newton polytope of polynomials of a chosen degree, reflected
in each orthant.  Due the fact that we do not use only polynomial and
therefore Newton Polytope, we did not manage to use \cite{Bra21,Sch91} to
shorten the proof.

\subsubsection{Decomposition methods using Descartes' rule of sign}

For univariate polynomials, there is a similar available criteria: Descartes'
rule of sign. It allows to ensure that a polynomial as at most one root in a
given interval. There are many attempts to generalise this rule to the
multivariate case. Most of these work
\cite{Itenberg1996,Lagarias1997,Bihan2017} only consider the zero dimensional
case. They give an upper bound of the number of solutions. If this upper bound
is one, it would allow to isolate each solution. A counter-example to Roy and
Itenberg conjecture \cite{Itenberg1996} given by Li and Wang
\cite{LiWang1998}, means that even for the zero dimensional case, there is no
known generalisation of Descartes' rule of sign that works in all cases.  In a
recent work \cite{feliu2022generalizing} Descates rules of sign is generalised
to arbitrary hypersurfaces, but only give a bound to the number of components
which is not enough to compute the topology of hypersufaces.

\subsubsection{Other decomposition methods}

Our approach is similar to many other algorithms that work by subdivising the
ambient space in hypercubes or simplices. See \cite{BPR06} for a book covering
most algorithms in the domain. In dimension 3, Marching cube algorithms usually
ensure correct topology between the piecewise trilinear function given by the
value at the vertex of each hypercube and the produced piecewise affine
variety \cite{Theisel2002, Grosso2016}. There are not many algorithms for
hypersurfaces in arbitrary dimension, we may cite \cite{Castelo2021}.

For works that provide exact algorithms using decomposition, in the case of polynomials
for 2D or 3D curves or 3D surfaces, we may cite \cite{Cha94,BS09}. Some of
these works, like \cite{LMP07} also use Bernstein basis in the
implementation. This latter also uses Descartes' rule of sign of the partial
derivative of the polynomials.

As we mentioned it previously, we think a key property of our criteria to stop
a decomposition is its global nature. Moreover, we are not aware of any
criteria that would work in arbitrary dimension and codimension.

\subsubsection{Algebraic methods}

Other algorithms, which are exact to compute topological invariants of algebraic
variety, are roughly based on the decidability of the theory of real closed
fields. Thus, they are more \emph{algebraic}, using cylindrical
decompositions, Groebner basis, resultants, \dots. They have the main
advantage of being able to deal with arbitrary singularities, but are of a
very different nature. Moreover, very few of these algorithms have free
implementations and we could not compare them with our algorithm, for instance
on the computation of Betti numbers of random cubic in dimension up to 6.

Some available implementations, that could fit in this \emph{algebraic} category, are limited to curves
and surfaces like \verb!bertini_real! \cite{BBHHSW17} or to the zero
dimensional case like \verb!msolve! \cite{Ber21}. Those two are probably the
best ones available. It is worth noticing that, unlike most available
implementation for curve and surfaces, \verb!bertini_real! does not limit the
number of variables and allow to compute a surface embedded in a space of high
dimension. However, unlike \verb!msolve!, \verb!bertini_real! is only using
arbitrary precision floating point arithmetic and is not an exact algorithm.

Currently the above cited algorithms support some management of
singularities and seem faster than our algorithm. This is to be expected as it
is natural that more general algorithms are slower than specialised
ones. Moreover, our search for a simplicial decomposition that satisfies our
criteria is in an early stage.  We outline below some directions of research
that could allow our criteria to be used in an algorithm that could compete with
the state of the art algorithms.

Remark: a lot of the algorithms found in the literature are not freely available
or hard to install. We only succeeded to install \verb!bertini_real! and \verb!msolve!!

\subsection{Further theoretical research}

Apart from proving or disproving the afford mentioned conjectures,
it would be nice to provide a complexity analysis for our algorithm. This will
require (as for some algorithms searching for roots of univariate polynomials), a
measure of regularity of the system of polynomials. This bound will probably
be very bad because it will assume the worst everywhere, but it would still
be interesting.

We would like to extend our work to product of projective spaces,
weighted projective spaces or compact of ${\mathbb R}^n$ with a border
condition allowing the variety to meet the border.

For singular varieties, it is likely that our criteria is still correct for
singularities which are affine subvariaties, provided that all singularities
of dimension $m$ are entirely covered by some faces of the simplicial
decomposition of dimension $m$. For instance, isolated points must be among
the vertices of the decomposition. Then, the only modification of our criteria
is to ignore the singular faces and it seems to work. We tested this on
isolated singularities, and this seems to work.  Singularities which are not
affine subvarieties seem must more challenging.

\subsection{Further implementation research}

The main problem with the current implementation is that we do not know yet
what are the best triangulations for our criteria, especially in the case of
codimension greater than one. For instance, if we impose the vertices of the
decomposition, what is the best triangulation to try to meet the criteria?
Currently we use the convex hull of the vertices projected on the unit
sphere. This is similar to Delauney's triangulation. This is frame independant,
but there is no reason that such triangulations are the best to meet our
criteria for a given polynomial system.  The same is true for the choice of
vertices. Currently, we favour critical points as it seems to give good
results, but this is not frame independant and does not always give enough
points and we don't really known what other points to choose.

We should also note that our implementation is written in OCaml using functors
to parameterise the implementation by the representation of numbers, because
it allows for rapid prototyping. A C implementation optimised for speed could
gain a factor 2 or 3 and parallelisation could allows to gain a factor 10 and
should be possible using OCaml 5.

We think it is possible to reach computing time matching those of existing
decomposition algorithms for curves and surfaces.

Another way to improve the efficiency is to combine our criteria with variables
elimination techniques. An idea would be to perform \emph{easy} eliminations,
before using our algorithm. For instance, one could eliminate a variable if it
occurs only in one monomial of some polynomials. The current implementation is
not even doing elimination of linear equations! But this is planed.

\subsection{Thanks}

We thank St\'ephane Simon for showing us his marching cube implementation,
25 years ago, starting our interest in this research topic. We thanks Ilia
Itenberg for several discussions, in particular about Viro's method. Finally,
we heartily thanks Frédéric Mangolte for the lengthy discussions on this
research, his comments and great help.


\section{Notation and convention}\label{convention}

Here are a few notation we use:

\begin{itemize}
\item $\hull{S}$ denotes the convex hull of a subset $S$ of $\mathbb{R}^n$.
\item When $f: \mathbb{R}^n \to \mathbb{R}$ is
  differentiable in all direction, we denote $\nabla(f)(x)(v)$ the
  differential of $f$ at $x$ in the direction $v$. In general, we only have
  $\nabla(f)(x)(\lambda v) = \lambda \nabla(f)(x)(v) $  for $\lambda > 0$ as
  $v \mapsto \nabla(f)(x)(v)$ may be non linear.
\item When $f: \mathbb{R}^n \to \mathbb{R}^m$ is differentiable
  $\nabla(f)(x)$ will denote the $m \times n$ Jacobian matrix.
\end{itemize}

\paragraph*{Simplicial decomposition}

In section \ref{compact}, we consider  simplicial decomposition of a compact polyhedron $\mathbb{K} \subset \mathbb{R}^n$. By
simplicial decomposition, we mean,  as in \cite{Ris93}, a family of simplices
$(S_i)_{i \in I}$ such that:
\begin{itemize}
\item $\mathbb{K} = \cup_{i \in I} S_i$
\item $\forall i,j \in I, i\neq j, S_i \cap S_j$ is a simplex of dimension at most $n-1$ which is the
  common face of $S_i$ and $S_j$.
\item $\forall i,j \in I, i\neq j, \mathring{S}_i \cap \mathring{S}_j = \emptyset$.
\end{itemize}

\paragraph*{Decomposition in simplicial cones}

\begin{defi}
  A \emph{simplicial cone}, is a set $S$, that is defined from a simplex $S'$
  that do not contain $0$ by
  $$S = \{ \lambda x, x \in S', \lambda > 0 \}$$
\end{defi}

In section \ref{projective} and after, we consider decomposition in simplicial
cones of $\mathbb{R}^{n+1}$. We mean a family of simplicial cone
$(S_i)_{i \in I}$ such that:
\begin{itemize}
\item $\mathbb{K} = \cup_{i \in I} S_i$
\item $\forall i,j \in I, i\neq j, S_i \cap S_j$ is a simplicial cone of dimension at most $n$ which is the
  common face of $S_i$ and $S_j$.
\item $\forall i,j \in I, i\neq j, \mathring{S}_i \cap \mathring{S}_j = \emptyset$.
\end{itemize}

\paragraph*{Bernstein basis}

In section \ref{implementation} we refer to Bernstein basis. In the case of
homogeneous polynomials of degree $d$, it is
 $$\left(\frac{d!}{\alpha!} x^\alpha\right)_{\alpha \in
    \mathbb N^d,  \Sigma_i  \alpha_i = d} \text{ where } x^\alpha = \Pi_i
  x_i^{\alpha_i} \text{ and } \alpha! = \Pi_i \alpha_i!$$

The key property of Bernstein basis it that its value in the unit simplex
lies in the convex hull of the coefficients. It is an immediate consequence of
De Casteljau algorithm to compute the value of the polynomial as a
barycenter. We also use this property with the gradient of a polynomial, seen
as a polynomial whose coefficients are vectors.


\section{Hypersurfaces on a compact polyhedron}\label{compact}

Here is a first theorem for an hypersurface which is enclosed in the interior
of a compact polyhedron of $\mathbb{R}^n$. This hypothesis seems essential and
unnatural, but will disappear when we consider the entire projective space of
dimension $n$.

\begin{theo}\label{thcompact}
  Let $\mathbb{K} \subset \mathbb{R}^n$ be a compact polyhedron.
  Let $(S_i)_{i\in I}$ be a simplicial decomposition of $\mathbb{K}$.
  Let $f : \mathbb{R}^n \to \mathbb{R}$ of class $C^1$.
  Let $V = \{ x \in \mathbb{K}, f(x) = 0\}$ be the zero locus of $f$
  restricted to $\mathbb{K}$.
  Assume that $V \cap \partial\mathbb{K} = \emptyset$ \reflabel{cond0}.

  We define $\tilde{f} : \mathbb{K} \to \mathbb{R}$ the piecewise affine
  function such that for all $i \in I$, $\restr{\tilde f}{S_i}$ is affine
  and for any $v$ vertex of $S_i$, we have
  $f(v) = \restr{\tilde f}{S_i}(v)$.
  We define the following:
  \begin{itemize}
  \item $\tilde V = \{ x \in \mathbb{K}, \tilde f(x) = 0\}$ the zero locus of
  $\tilde f$.
  \item $\mathbb{K}(f) = \{ x \in \mathbb{K}, f(x)\tilde{f}(x) \leq 0\}$.
  \item $\tilde \nabla f(x) = \{ \nabla \restr{\tilde f}{S_i}(x), x \in S_i
    \} \subset \mathbb{R}^n$
  \item $G(f,x) = \{ \nabla f(x) \} \cup \tilde \nabla f(x) \subset \mathbb{R}^n$
  \end{itemize}
  If the condition below holds, then $V$ and $\tilde V$ are isotopic:
  \begin{eqnarray}
    \forall x \in \mathbb{K}(f), 0 \notin \hull{G(f,x)} \label{eq1}
  \end{eqnarray}
\end{theo}

\begin{proof}
  Assume the definitions and hypotheses of the theorem.
  For any $x \in \mathbb{R}^n$ and $\epsilon >
  0$, we define
  $$ G(f,x,\epsilon) = \bigcup_{y \in \mathbb{K}(f), \|y - x\| < \epsilon} G(f,y) $$

  Remark: we need to define $G(f,x,\epsilon)$ for $x \in \mathbb{R}^n$ because
  the convolution product below will cover the border of
  $\mathbb{K}(f)$. Clearly, for points too far from $\mathbb{K}(f)$, we have
  $G(f,x,\epsilon) = 0$.

  We now prove that there exists $\epsilon > 0$ such that
  \begin{eqnarray}
      \forall x \in \mathbb{K}(f), 0 \notin \hull{G(f,x,\epsilon)} \label{eq2}
  \end{eqnarray}
  We proceed by contradiction and choose a sequence $(x_n)_{n \in \mathbb{N}}$ in
  $\mathbb{K}(f)$ such that $0 \in \hull{G(f,x_n,\frac{1}{n})}$. As
  $\mathbb{K}(f)$ is compact, we can assume that $x_n$ converges to $x_\infty \in
  \mathbb{K}(f)$. Let us define
  $$\delta = \min_{i\in I, x_\infty \not\in S_i} \dist(x_\infty,S_i)$$
  We have
  for any $y \in \mathbb{K}$,
  $\|y - x_\infty\| < \delta$ and $y \in S_i$ implies $x_\infty \in S_i$ and therefore
  $\tilde \nabla f(y) \subset \tilde \nabla f(x_\infty)$.
  Thus, for $\delta n > 1$, we have
  $$0 \in \hull{\{ \nabla f(y), y \in \mathbb{K}, \|x_n - y\| < \frac{1}{n} \}
    \cup \tilde \nabla f(x_\infty)}$$

  The set $\{ \nabla f(y), y \in \mathbb{K}, \|x_n - y\| < \frac{1}{n} \}$
  converges to the singleton $\{\nabla f(x_\infty)\}$ for the Haussdorf metric
  and $\mathcal H$ is continuous for that metric. Hence,
  $$\hull{\{ \nabla
    f(y), y \in \mathbb{K}, \|x_n - y\| < \frac{1}{n} \} \cup \tilde \nabla
    f(x_\infty)} \longrightarrow \hull{G(f,x_\infty)} \text{ when } n \longrightarrow +\infty$$
  This implies $0 \in G(f,x_\infty)$, because it is a closed set,
  which contradicts (\ref{eq1}).

  By the geometric form of Hahn-Banach, we
  can find $N: \mathbb{R}^n \to \mathbb{R}^n$ such that
  $$\forall x \in \mathbb{R}^n, \forall v \in G(f,x,\epsilon), N(x).v > 0$$

  Let us choose a function $\mu : \mathbb{R}^n \to \mathbb{R}_+$ of $C^\infty$
  class, with support in the sphere of radius
  $\epsilon$ and such that $\int_{\mathbb{R}^n} \mu(u) \mathrm du = 1$. We
  define:

  $$N'(x) = N \star \mu = \int_{\mathbb{R}^n} N(u - x) \mu(u) \mathrm du$$

  $N'$ is of $C^\infty$ class on $\mathbb{R}^n$. Let us consider $v \in
  G(f,x)$ for $x \in \mathbb{K}$, if $\|u\| < \epsilon$ and therefore $\|x - (u - x)\| < \epsilon$,
  we have $v \in G(f,u-x,\epsilon)$ which implies $N(u - x).v > 0$. This
  establishes:

  \begin{eqnarray}
    \forall x \in \mathbb{K}(f), \forall v \in G(f,x), N'(x) . v > 0 \label{eq3}
  \end{eqnarray}

  The next step is to consider the maximal integral curves of $N'$. By
  The Cauchy-Lindel\"of-Lipshitz-Picard theorem those curves exist, are
  unique in $\mathbb{K}(f)$ and continuous in the initial conditions.

  For $t \in [0,1]$ we define $f_t : \mathbb{K} \to \mathbb{R}$, such that
  $f_t(x)  = t f(x) + (1 - t) \tilde f(x)$. We remark that
  $\tilde f$ is differentiable in every direction and that
  the differential of $\tilde f$ at $x$ in the direction $D$ is given
  by $D.V$ for some $V \in \tilde \nabla f(x)$. Remark the differential in the
  direction $D$ and $-D$ may be different.

  Therefore, the differential of $f_t$ at $x \in \mathbb{K}(f)$ in the
  direction $N'(x)$ is given by the expression $N'(x).(t \nabla f(x) + (1 - t)
  V)$ for some $V \in \tilde \nabla f(x)$ and is therefore positive by
  (\ref{eq3}).

  This means that the functions $f_t(x)$ are increasing along an integral
  curve of $N'$ and therefore each integral curve meet the variety
  $V_t = \{ x \in \mathbb{R}, f_t(x) = 0\}$ for $t \in [0,1]$ in at most one
  point. Remark that $V_t \subset \mathbb{K}(f)$.

  To finish the proof we must show that the maximal integral curves of $N'$ have their
  extremity in the border of $\mathbb{K}(f)$, which are points $x$ with either
  $f(x) = 0$ or $\tilde f(x) = 0$, by the condition ($\ref{cond0}$). This is
  true because we can find $K > 0$ such that $\forall x \in \mathbb{K}(f), K
  < N'(x)$ by compacity and regularity of $N'$. This means that a maximal integral
  curve of $N'$ will join the border of $\mathbb{K}(f)$ on an interval
  $[t_1,t_2]$ for $t_2 - t_1 < \frac{M}{K}$ where $M$ is an upper bound of
  both $f$
  and $\tilde{f}$ on $\mathbb{K}$.

  Therefore, $(x,t) \mapsto f_t(x)$ is the wanted isotopy.
\end{proof}

\section{Hypersurfaces on the projective space}\label{projective}

We now give a condition to establish that an hypersurface in $\mathbb{S}^n$
the unit sphere of $\mathbb{R}^{n+1}$
defined by a positively homogeneous $C^1$ function in $n+1$ variables is istopic to a
variety defined by a piecewise linear function on
$\mathbb{R}^{n+1}$. We state the theorem on the unit sphere $\mathbb{S}^n$
because it is simpler to write the condition than working on the projective space.

\begin{theo}\label{thprojective}
  Let $(S_i)_{i\in I}$ be a decomposition of $\mathbb{R}^{n+1}$ in simplicial
  cones with vertices on $\mathbb{S}^n$, the
  unit sphere of $\mathbb{R}^{n+1}$.
  Let $p : \mathbb{R}^{n+1} \to \mathbb{R}$ be a positively homogeneous $C^1$
  function of degree $d$.
  (i.e. $p(\lambda x) =
  \lambda^{d} p (x)$ for any $\lambda \in \mathbb{R}_+$ and $x \in
  \mathbb{R^{n+1}}$).
  Let $V = \{ x \in \mathbb{S}^n, p(x) = 0\}$ be the zero locus of $p$
  restricted to $\mathbb{S}^n$.

  We define $\tilde{p} : \mathbb{R}^{n+1} \to \mathbb{R}$ the piecewise linear
  function such that for all $i \in I$, $\restr{\tilde p}{S_i}$ is linear
  and for any $v$ vertex of $S_i$, we have
  $p(v) = \restr{\tilde p}{S_i}(v)$.
  We define the following:
  \begin{itemize}
  \item $\tilde V = \{ x \in \mathbb{S}^n, \tilde p(x) = 0\}$ the zero locus of
  $\tilde p$.
  \item $\mathbb{K}(p) = \{ x \in \mathbb{S}^n, p(x)\tilde{p}(x) \leq 0\}$.
  \item $\tilde \nabla p(x) = \{ \nabla \restr{\tilde p}{S_i}(x), x \in S_i
    \} \subset \mathbb{R}^{n+1}$
  \item $G(p,x) = \{ \nabla p(x) \} \cup \tilde \nabla p(x) \subset \mathbb{R}^{n+1}$
  \end{itemize}
  If the following condition (\ref{peq1}) holds, then $V$ and $\tilde V$ are isotopic:
  \begin{eqnarray}
    \forall x \in \mathbb{K}(p), 0 \notin \hull{G(p,x)} \label{peq1}
  \end{eqnarray}
  This implies that the projective varieties associated to $V$ and $\tilde V$
  are isotopic if the simplicial decomposition is stable by the symmetry $x
  \mapsto -x$.
\end{theo}

\begin{proof}
  By taking $\mathbb{K} = \mathbb{S}^n$, we can use exactly the same
  definitions and reasoning as the proof of theorem \ref{thcompact} until
  the definition of $N'$ of $C^\infty$ class satisfying

  \begin{eqnarray}
    \forall x \in \mathbb{K}(p), \forall v \in G(p,x), N'(x) . v > 0 \label{peq3}
  \end{eqnarray}

  We must change the definition of $p_t(x)$, using  $d$ the degree of $p$:
  \begin{eqnarray*}
    \overline{p}(x) &=& \tilde p(x \|x\|^{d-1}) \cr
    p_t(x) &=& t p(x) + (1 - t)\overline p(x)\cr
  \end{eqnarray*}

  We have $p_t(\lambda x) = \lambda^d p_t(x)$ for all $\lambda \in
  \mathbb{R}_+$.
  As in the previous proof, $\tilde p(x)$ has a differential in any direction
  $v$. Hence, the functions $\overline p$ and $p_t$ are differentiable in
  any direction, hence they satisfy Euler relation:
  $$\nabla \overline p(x)(x) = d  p(x)  \hspace{3cm} \nabla p_t(x)(x) = d  p(x)$$

  We need more precision for the derivative of $\tilde p$: if $x,v \in
  \mathbb{R}^{n+1}$ we can define $i(x,v) \in I$ such that $x \in S_{i(x,v)}$
  and $x + hv \in S_{i(x,v)}$ for all $h > 0$ small enough. We also define
  $S(x,v) = S_{i(x,v)}$. Using this notation, the gradient of $\tilde
  p(x)$ in the direction $v$ is given by
  \begin{eqnarray}
  \nabla \tilde p(x)(x) &=& \nabla \restr{\tilde
    p}{S(x,v)}(x).v
  \end{eqnarray}
  But as $\restr{\tilde p}{S(x,v)}$ is linear, its gradient
  is constant and we can simply write $\nabla \restr{\tilde
    p}{S(x,v)}.v$. Remark: $i(x,v)$ is not uniquely defined, but the
  choice of index in $I$ does not change the value of the differential in a
  given direction.

  Using these notations, we compute:
  \begin{eqnarray*}
    \nabla \overline{p}(x)(v) &=& d \nabla \restr{\tilde
      p}{S(x,v)}.v \|x\|^{d-1} \cr
    \nabla p_t(x)(v) &=& (t \nabla p(x) +
      (1-t) d \nabla \restr{\tilde p}{S(x,v)} \|x\|^{d-1}).v
  \end{eqnarray*}

  We now prove that for $x \in \mathbb{K}(p)$, $N'(x)$ is not normal to
  the unit sphere at $x$. Let us choose $x \in \mathbb{K}(p)$, we can
  find $t \in [0,1]$ such that $p_t(x) = t p(x) + (1 - t)\overline p(x) = 0$
  (take $t = 0$ if $p(x) = \tilde p(x) = 0$ and $t = \frac{\overline
    p(x)}{\overline p(x) - p(x)}$ otherwise). This is well defined because in
  $\mathbb{K}(p)$ we can only have $p(x) = \overline p(x)$ if $p(x) = \tilde
  p(x) = 0$. Let us define
  $$V = t \nabla p(x) +  (1-t) d \nabla \restr{\tilde p}{S(x,v)}$$
  The gradient $p_t(x)$ in the direction of $x$,
  is given by $V.x = p_t(x) = 0$ by Euler relation.
  $V \in (t + (1-t)d) G(p,x)$ and $t + (1-t)d = d - t(d-1) > 0$ for $t \in
    [0,1]$.
    This means that $N'(x)$ can not be normal to
  $\mathbb{S}^n$ at $x$ as this would imply $V.N'(x) = 0$ which is impossible
  by (\ref{peq3}). This ends the proof that  $N'(x)$ is not normal to
  the unit sphere at $x$.

  We can now define for all $x \in \mathbb{K}(p)$, $N^T(x) = N'(x) -
  (N'(x).x)x$, the projection of $N'(x)$ on the hyperplane tangent to
  $\mathbb{S}^n$ at $x$. For $x \in \mathbb{K}(p)$,
  as the polyhedra $S_i$ are simplicial cones we can assume $S(x,N'(x)) =
  S(x,N^T(x))$ that we will simply write $S(x)$. Indeed, for $h > 0$ small
  enough, $x + hN'(x)$ and $x + hN^T(x)$ belong to the same simplicial cone $S_i$
  because they only differ by a vector in the direction of $x$ and $S_i$ is a cone.

  Using this notation, for a point $x \in \mathbb{K}(p)$ such that $p_t(x) =
  0$, the gradient of $p_t(x)$ in direction $N'(x)$
  and $N^T(x)$ verify:
  \begin{eqnarray}
    \nabla p_t(x)(N'(x)) &=& (t \nabla p(x) + (1-t) d \nabla \restr{\tilde
      p}{S(x)} ).N'(x) \cr
    &>& 0 \hfill \text { by (\ref{peq3}) } \cr
    \nabla p_t(x)(N^T(x)) &=& (t \nabla p(x) + (1-t) d \nabla \restr{\tilde
      p}{S(x)} ).N^T(x) \cr
    &=& (t \nabla p(x) + (1-t) d \nabla \restr{\tilde
      p}{S(x)} ).(N'(x) - (x.N'(x))x) \cr
    &=& (t \nabla p(x) + (1-t) d \nabla \restr{\tilde
      p}{S(x)} ).N'(x) \hfill \text{ by Euler relation}\cr
    &>& 0 \label{peq4}
  \end{eqnarray}

  Let $\gamma : J \to \mathbb{K}(p)$ be a maximal integral curve of
  $N^T$. This means $\gamma'(u) = N^T(\gamma(u))$.
  There are particular cases where
  $\gamma$ is reduced to one point $x$ when $p(x) = \tilde p(x) = 0$.

  In all other cases, as $\tilde p$ is derivable in all directions, $u \mapsto
  \tilde p(\gamma(u))$ is derivable (but not necessarily of $C^1$
  class). Similarly $u \mapsto \overline{p}(\gamma(u))$ and $u \mapsto
  p_t(\gamma(u))$ for $t \in [0,1]$ and $t(u) = \frac{\overline p(u)}{\overline p(u)
    - p(u)}$ are derivable. Moreover, a point $x \in \mathbb{K}(p)$ with $p(x) = 0$ or $\tilde p(x) = 0$ can only be
  at the extremity of $\gamma$, otherwise, by (\ref{peq4}), $\gamma$ would leave
  $\mathbb{K}(p)$. This means that $p$ and $\tilde p$ have constant
  sign on $\gamma$ and may be null only on the extremity. This implies that
  $\overline{p}(\gamma(u)) - p(\gamma(u))$ is of contant signe along such a
  curve $\gamma$.

  We have:
  \begin{eqnarray}
    p_{t(u)}(\gamma(u)) &=& 0\cr
    (p(\gamma(u)) - \overline{p}(\gamma(u))) t'(u) + \nabla
    p_{t(u)}(\gamma(u))(\gamma'(u)) &=& 0\cr
    (\overline{p}(\gamma(u)) - p(\gamma(u))) t'(u) &=&  \nabla
    p_{t(u)}(\gamma(u))(N^T(\gamma(u))) \cr
    (\overline{p}(\gamma(u)) - p(\gamma(u)))t'(u) &>& 0
  \end{eqnarray}

  This means that $t(u)$ is monotonous along any curve $\gamma$ that is not
  reduced to one point. As in the previous proof, $N^T(\gamma(u))$ is never
  null and is minored by some constant $K > 0$, thus the extremity of maximal integral
  curve $\gamma$ will necessarily be a point where $p$ is null and another
  point where $\tilde p$ is null.

  This means that $p_t$ defines the isotopy we are looking for.
\end{proof}


\section{Increasing codimension}\label{codimension}

We propose the following conjecture for several positively homogeneous $C^1$ funciton:

\begin{conj}\label{conj}
  Let $(S_i)_{i\in I}$ be a decomposition of $\mathbb{R}^{n+1}$ in simplicial
  cones with vertices on $\mathbb{S}^n$, the
  unit sphere of $\mathbb{R}^{n+1}$.
  Let $p = (p_1,\dots,p_m)$ be a family of $m \leq n$ positively homogeneous
  $C^1$ functionx from $\mathbb{R}^{n+1}$ to $\mathbb{R^n}$ of respective
  degree $(d_1,\dots,d_m)$ not necessarily equal (i.e. $p_i(\lambda x) =
  \lambda^{d_i} p_i (x)$ for any $\lambda \in \mathbb{R}_+$ and $x \in \mathbb{R^{n+1}}$).
  Let $V = \{ x \in \mathbb{S}^n, p(x) = 0\}$ be the zero locus of $p$
  restricted to $\mathbb{S}^n$.

  We define $\tilde{p} : \mathbb{R}^{n+1} \to \mathbb{R}^m$ the piecewise linear
  function such that for all $i \in I$, $\restr{\tilde p}{S_i}$ is linear
  and for any $v$ vertex of $S_i$, we have
  $p(v) = \restr{\tilde p}{S_i}(v)$.
  We define the following:
  \begin{itemize}
  \item $\tilde V = \{ x \in \mathbb{S}^n, \tilde p(x) = 0\}$ the zero locus of
  $\tilde p$.
  \item $\mathbb{K}(p) = \{ x \in \mathbb{S}^n, \forall i \in \{1,\dots,m\}, p_i(x)\tilde{p}_i(x) \leq 0\}$.
  \item $\tilde \nabla p(x) = \{ \nabla \restr{\tilde p}{S_i}(x), x \in S_i
    \} \subset \mathrm{Mat}_{m,n+1}(\mathbb{R})$
  \item $G(p,x) = \{ \nabla p(x) \} \cup \tilde \nabla p(x)  \subset \mathrm{Mat}_{m,n+1}(\mathbb{R})$
  \end{itemize}
  If the condition below holds, then $V$ and $\tilde V$ are isotopic:
  \begin{eqnarray}
    \forall x \in \mathbb{K}(p), \forall A \in \hull{G(p,x)}, A \text{ is of
      maximal rank} \label{ceq1}
  \end{eqnarray}
\end{conj}

This conjecture is the natural generalisation of theorem \ref{thprojective}
and our implementation described in the next section suggest that it might be true.

Unfortunately, to adapt the proof of the previous section, we need a result
analogous to Hahn-Banach for convex set of full rank matrices. This would give
the countepart of the vector $N(x)$ in the previous proof and also the
certificate we need for the implementation.

Here is the expected result that seems unknown and that we could not prove
neither disprove:

\begin{conj}\label{matconj}
  Let $1 < m \leq n$ two natural numbers, let $S \subset
  \mathrm{Mat}_{m,n}(\mathbb R)$ a
  convex set of matrices of rank $m$. There exists a matrix $M \in
  \mathrm{Mat}_{m,n}(\mathbb R)$ such that $M \tr A + A \tr M$ is symmetric definite and
  positive for all $A \in S$.
\end{conj}

We could not prove that conjecture \ref{matconj} implies conjecture \ref{conj}
but we feel that it is a key element of the proof. A way to prove this
implication would be to ensure that $N^T : \mathbb{K}(p) \to
\mathrm{Mat}_{m,n+1}(\mathbb{R})$ satisfies the Schwartz condition.
This means we should have that the
derivative of $N_i^T$ in the direction $N_j^T$ should be equal to the
derivative of $N_j^T$ in the direction $N_i^T$.
Then, we could construct unique integral hypersurfaces of $N^T$ and
probably finish the proof. To to this, an idea if to build the kernel used in the
convolution product defining $N'$ by solving a partial differential equation...

As we can not prove constructively conjecture \ref{matconj}, we have no
algorithm to test condition \ref{ceq1} that would give a certificate.
Therefore, for our implementation, we use a stronger condition using this definition:

\begin{defi}
  Let $1 < m \leq n$ two natural numbers, let $S \subset
  \mathrm{Mat}_{m,n}(\mathbb R)$ a
  set of matrices. $S$ is said to be
  \emph{strongly full rank} if
  $$\forall \sigma \in \{-1,1\}^m, 0 \notin \hull{\{\sigma A, A \in S\}}$$
\end{defi}

\begin{prop}
If $S \subset \mathrm{Mat}_{m,n}(\mathbb R)$ is strongly full rank, then
$\hull{S}$ contains only full rank matrices.
\end{prop}

\begin{proof}
  If $A \in \hull{S}$ is not full rank, then there exists
  $v \in \mathbb{R}^m$ such that $v A = 0$.
  Take $\sigma = (\sigma_1, \dots, \sigma_n)$ such that $\sigma_i = 1$ if $v_i
  \geq 0$ and $\sigma_i = -1$ otherwise and we find
  $0 \in \hull{\{\sigma A, A \in S\}}$.
\end{proof}

This stronger condition gives a weaker conjecture that correspond to our implementation:
\begin{conj}\label{wconj}
  Let $(S_i)_{i\in I}$ be a decomposition of $\mathbb{R}^{n+1}$ in simplicial
  cones with vertices on $\mathbb{S}^n$, the
  unit sphere of $\mathbb{R}^{n+1}$.
  Let $p = (p_1,\dots,p_m)$ be a family of $m \leq n$ positively homogeneous
  $C^1$ functionx from $\mathbb{R}^{n+1}$ to $\mathbb{R^n}$ of respective
  degree $(d_1,\dots,d_m)$ not necessarily equal (i.e. $p_i(\lambda x) =
  \lambda^{d_i} p_i (x)$ for any $\lambda \in \mathbb{R}_+$ and $x \in \mathbb{R^{n+1}}$).
  Let $V = \{ x \in \mathbb{S}^n, p(x) = 0\}$ be the zero locus of $p$
  restricted to $\mathbb{S}^n$.

  We define $\tilde{p} : \mathbb{R}^{n+1} \to \mathbb{R}^m$ the piecewise linear
  function such that for all $i \in I$, $\restr{\tilde p}{S_i}$ is linear
  and for any $v$ vertex of $S_i$, we have
  $p(v) = \restr{\tilde p}{S_i}(v)$.
  We define the following:
  \begin{itemize}
  \item $\tilde V = \{ x \in \mathbb{S}^n, \tilde p(x) = 0\}$ the zero locus of
  $\tilde p$.
  \item $\mathbb{K}(p) = \{ x \in \mathbb{S}^n, \forall i \in \{1,\dots,m\}, p_i(x)\tilde{p}_i(x) \leq 0\}$.
  \item $\tilde \nabla p(x) = \{ \nabla \restr{\tilde p}{S_i}(x), x \in S_i
    \} \subset \mathrm{Mat}_{m,n+1}(\mathbb{R})$
  \item $G(p,x) = \{ \nabla p(x) \} \cup \tilde \nabla p(x)  \subset \mathrm{Mat}_{m,n+1}(\mathbb{R})$
  \end{itemize}
  If $\forall x \in \mathbb{K}(p)$, $G(p,x)$ is strongly full rank, then $V$ and $\tilde V$ are isotopic.
\end{conj}


\section{Implementation}\label{implementation}

Obtaining an implementation from theorem \ref{thprojective} or conjecture
\ref{wconj} is not very difficult.

Let $(S_i)_{i\in I}$ be a decomposition of $\mathbb{R}^{n+1}$ in simplicial
cones with vertices on $\mathbb{S}^n$, the unit sphere of $\mathbb{R}^{n+1}$.
Let $p = (p_1,\dots,p_m)$ be a family of $m \leq n$ homogeneous polynomials
with $n+1$ variables.

Let us consider now a simplicial cone $F$ which is a subset of a face of
dimension $d_F$ of one
of the simplex $S_i$ for $i \in I$. $F$ could be reduced to a vertex when $d_F
= 0$ or could
be of dimension $d_F=n+1$. By doing a change of coordinates sending
the unit simplex of dimension $d_F$ to $F$ and writing the resulting polynomials in the Berstein
bases, we can use the fact that the value of polynomials or their
differentials are in the convex hull of the coefficients to check that the
condition of our theorem \ref{thprojective} or conjecture
\ref{wconj} holds in $F$.

This leads to the following procedure:

\begin{algo}[{\tt TEST\_FACE}] \;

  \noindent Inputs:
  \begin{itemize}
  \item $(S_i)_{i\in I}$ a decomposition of $\mathbb{R}^{n+1}$ in simplicial
    cones with vertices on $\mathbb{S}^n$.
  \item $p = (p_1,\dots,p_m)$ a family of $m \leq n$
    homogeneous polynomials with $n+1$ variables.
  \item $\tilde p = (\tilde p_1,\dots,\tilde p_m)$ the piecewise linear
    functions associated to $p$ and $(S_i)_{i\in I}$.
  \item a simplex $F$ of dimension $d_F$ which is included in a face of
    one of the $S_i$.
  \item a minimal size for simplices.
  \item a heuristic to split a simplex in $2$ (that may use all other inputs).
  \end{itemize}
  Algorithm:
  \begin{enumerate}
  \item Build the matrix $P$ sending the unit simplex of dimension $d_F$ to $F$
  \item Write $p(M(x))$ and $\tilde p(M(x))$ in the Berstein basis, this gives
    two families of homogeneous polynomials $q(x)$ and $\tilde q(x)$ with $d_F + 1$
    variables. Remark: if $F$ is
    a vertex, it is equivalent to evaluating the polynomials!
  \item If there is $1 \leq i \leq m$ such that all coefficients of $q_i$ and
    $\tilde q_i$ have the same sign, return \verb#TRUE# because $F$ does not meet $\mathbb{K}(p)$.
  \item Otherwise, compute the list $L$ such that $\{S_l, l \in L\}$ is the set all
    simplicies that contains $F$.
  \item Write $\nabla p(M(x))$ and, for all $l \in L$, $\nabla \restr{\tilde
    p}{S_l}(M(x))$ in the Berstein bases. Define the set $A$ of $m \times (n + 1)$
    matrices which are the coefficients of those polynomials. If for all
    $\sigma \in \{-1,1\}^m$ we have $0 \notin \hull{\{\sigma M, M \in A\}}$
    return \verb#TRUE# because $A$ is strongly full rank.
  \item Otherwise, if $F$ is not too small, subdivide $F$ in $F_1$ and $F_2$ and
    recursively call the procedure \verb#TEST_FACE# on $F_1$ and $F_2$ and
    return \verb#TRUE# if both calls return \verb#TRUE#.
  \item Otherwise, if $F$ was too small, return \verb#FALSE#.
  \end{enumerate}
\end{algo}

\noindent Using this procedure, we can implement our main loop:

\begin{algo}[{\tt MAIN\_LOOP}] \;

  \noindent Inputs:
  \begin{itemize}
  \item $(S_i)_{i\in I}$ a decomposition of $\mathbb{R}^{n+1}$ in simplicial
    cones with vertices on $\mathbb{S}^n$.
  \item $p = (p_1,\dots,p_m)$ a family of $m \leq n$
    homogeneous polynomials with $n+1$ variables.
  \item $\tilde p = (\tilde p_1,\dots,\tilde p_m)$ the piecewise linear
    functions associated to $p$ and $(S_i)_{i\in I}$.
  \item a heuristic to refine the decomposition (that may use all other inputs).
  \item a minimal size for simplices.
  \item a heuristic to split a simplex in $2$
  \end{itemize}
  Algorithm:
  \begin{enumerate}
    \item For each face $F$ of one of the simplex call the procedure
      \verb#TEST_FACE#. If all calls return \verb#TRUE#, return $\tilde p$.
    \item If the procedure \verb#TEST_FACE# returned \verb#FALSE# on $F$,
      try to refine the decomposition, preferably in a way that splits $F$.
    \item Update $\tilde p$ to the new decomposition.
    \item Call back \verb#MAIN_LOOP# with the refined subdivition and new
      $\tilde p$.
  \end{enumerate}
\end{algo}

\noindent Here is the entry point of our implementation:

\begin{algo}[{\tt MAIN}]  \;

  \noindent Inputs:
  \begin{itemize}
  \item $p = (p_1,\dots,p_m)$ a family of $m \leq n$
    homogeneous polynomials with $n+1$ variables.
  \item a heuristic to refine a decomposition of $\mathbb{R}^{n+1}$ in simplicial cones.
  \item a minimal size for simplices.
  \item a heuristic to split a simplex in $2$.
  \end{itemize}
  Algorithm:
  \begin{enumerate}
  \item Build $(S_i)_{i\in I}$ a decomposition of $\mathbb{R}^{n+1}$ in simplicial
    cones with vertices on $\mathbb{S}^n$.
  \item Build $\tilde p = (\tilde p_1,\dots,\tilde p_m)$ the piecewise linear
    functions associated to $p$ and $(S_i)_{i\in I}$.
  \item Call the procedure \verb#MAIN_LOOP#.
  \item If it return $\tilde p$, build the piecewise affine projective variety of
    equation $\tilde p(x) = 0$ and return it.
  \end{enumerate}
\end{algo}

\begin{prop}
  Given $p = (p_1,\dots,p_m)$ a family of $m \leq n$
    homogeneous polynomials with $n+1$ variables,
  if there is only one polynomial or if conjecture \ref{wconj} is true, the
  above algorithm loops or returns a piecewise affine projective variety that
  is isoptopic to the variety defined by $p(x) = 0$.
\end{prop}

\begin{proof}
  This is a consequence of the definitions and the properties of Bernstein
  basis. It is important to note that the property ``$0$ is the convex hull''
  used in the algorithm is invariant by a linear change of variable.
\end{proof}

\begin{prop}
   Let $p = (p_1,\dots,p_m)$ a family of $m \leq n$ homogeneous polynomials
   with $n+1$ variables.  Assume that the matrix $\nabla p(x)$ is full rank for
   all $x \in \mathbb{S}^n$ such that $p(x) = 0$. Then, our algorithm
   terminates, if the heuristic to refine decompositions, when repeated, gives
   decompositions such that all simplicial cones have diameter that converges to $0$,
   when intersected with the unit sphere.
\end{prop}

\begin{proof}
  If the diameter of all simplicial cones restricted to $\mathbb{S}^n$ are
  small enough, then $p$ will be almost linear on each of them and the points
  where $p(x) = 0$ will be separated from the points where the matrix $\nabla
  p(x)$ is not full rank. This means that one of the two tests will always
  succeed in the procedure \verb#TEST_FACE#.
\end{proof}

\paragraph*{Test for $0$ in the convex hull}

Clearly we do not want to compute the convex hull to check for one point.
This problem is traditionally implemented as a reduction to linear
programming. We chose to implement it directly:

For a finite set of vector $A$, we try to minimise $\|N\|^2$ for $N = \sum_{V
  \in A} \alpha_V V$ with $\alpha_V > 0$ for all $V \in A$ and $\sum_{V \in A}
\alpha_V = 1$. We perform this minimisation by alterning two kinds of steps:
\begin{itemize}
  \item Linear steps: we solve a linear system to find a direction which
    is not always a direction of descent but that often offers rapid progress.
    This kind of steps may set some of the $\alpha_V$ to zero.
  \item Descent in the direction $V \in A$ if $N.V \leq 0$. It is easy to show that
    $\|\frac{N + \alpha V}{1 + \alpha}\|^2 < \|N\|^2$ in this case. This kind
    of steps increase $\alpha_V$, even if $\alpha_V = 0$. We stop if there is
    no such vector $V$ and we know that $0 \notin \hull{A}$.
  \item We stop if $\|v\|^2$ is too small (meaning we can probably reach $0$).
  \item It is important to avoid setting $\alpha_V
    = 0$ in a linear step, followed by a descent in the direction $V$. This
    yields to very slow progress. Thus, we do not select
    $V$ for descent if $\alpha_V$ was set to $0$ by the previous linear step.
\end{itemize}

This relatively simple algorithm works very well for this specific case and
might be the object of a separate publication in the near future. It is worth
noticing that our algorithm is not an interior point method nor a method that
stay on the border of the convex hull.

We mentioned the algorithm to make it clear that when we fail to find a
descent direction, we have $N.V > 0$ for all $V \in A$ and therefore $N$ is
the vector given by the geometric form of Hahn Banach theorem and considered
by the proof.

\paragraph*{Certificate and exact algorithm}

The above algorithm is implemented using 64 bits floating point numbers.
However, the procedure \verb#TEST_FACE# keeps a trace of the subdivision it
did as a binary tree and it also keeps in the leaf of the tree a boolean
giving the reason of
success: either \verb#TRUE# if the sign of the polynomials was constant or
\verb#FALSE# if the test for the convex
hull succeeded. In the latter case it also keeps the vectors $N$ given by the
algorithm for each value of $\sigma \in \{-1,1\}^m$.

Such a tree is associated to each face of the simplicial decomposition and
form a \emph{certificate}.

This allows to recheck the criteria using exact rational arithmetic and the
only operations are:
\begin{itemize}
\item change of coordinates in the polynomials,
\item scalar products and
\item comparisons.
\end{itemize}

As a result the final check of this certificate is fast (in practice faster
that the initial computation) and ensures an exact result (if conjecture
\ref{wconj} is true when codimension is greater than 1).


\section{Experiments}\label{experiments}

Note: all figures in this article use a projection of the projective space into
a sphere, so we see the entire variety.

\paragraph*{Study near singularity}

Our first example is with the family of quartic polynomials:

$$p_\epsilon(x,y,t) =   (x^2+y^2-t^2)^2 + \epsilon x y (x-y) (x+y)$$

When $\epsilon > 0$, the curve $p_\epsilon(x,y,t) = 0$ has four components and it converges
to a circle of double points when $\epsilon \rightarrow 0$. However, a simplicial decomposition with
only 24 triangles seems sufficient for any $\epsilon>0$. Only the number of
subdivision of each triangle increases when  $\epsilon \rightarrow 0$. Figure \ref{exp-quartic} is
the simplicial decomposition (in green) and the curve $\overline p_\epsilon(x,y,t) = 0$ we
get for some $\epsilon>0$ (in black):

\begin{figure}
  \begin{center}
    \includegraphics[width=0.4\linewidth,bb=0 0 1044 1039]{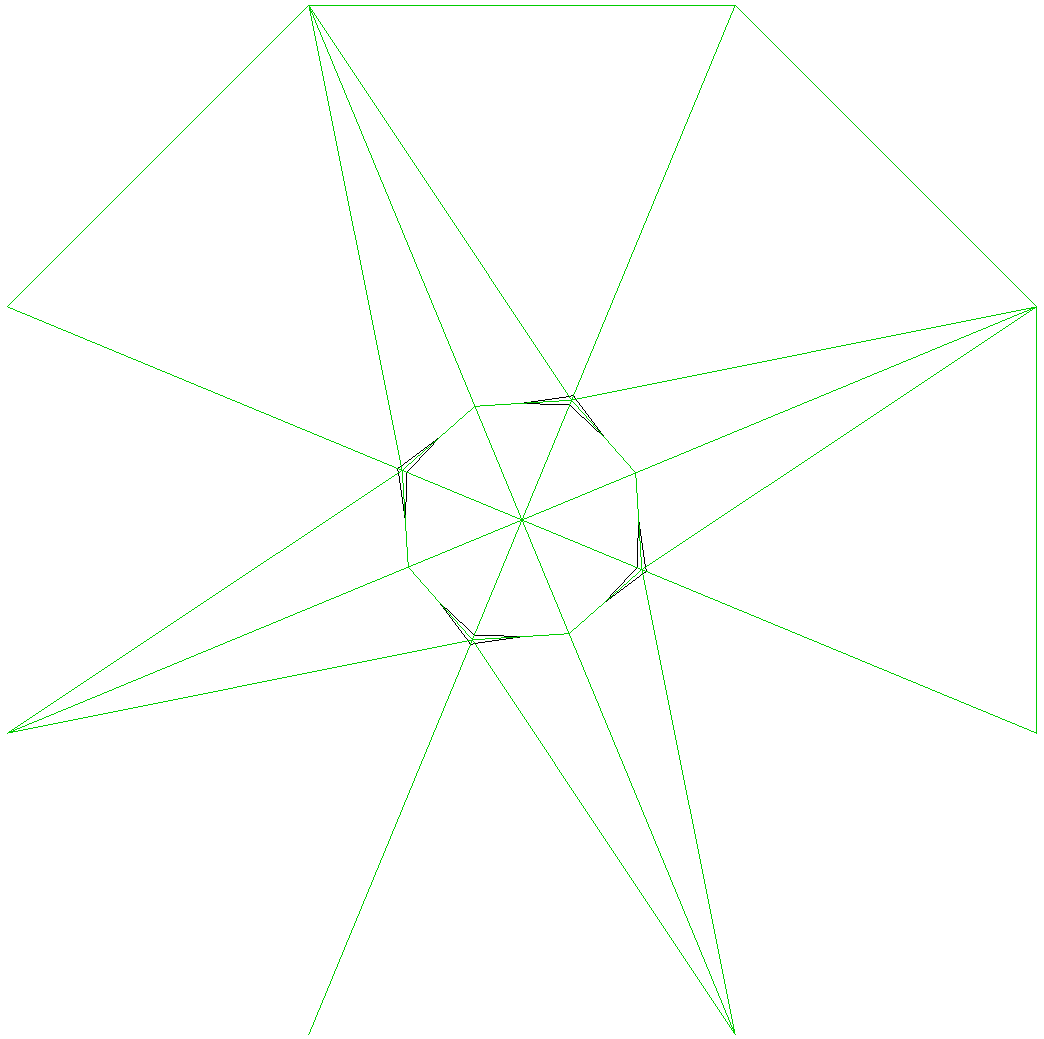}
    \caption{simplicial decomposition for $p_\epsilon(x,y,t) =   (x^2+y^2-t^2)^2 + \epsilon x y (x-y) (x+y)$}\label{exp-quartic}
  \end{center}
\end{figure}

We now give a table that gives for some value of $\epsilon$, the total computing time,
the time for the exact test using rational arithmetic and the
maximum number of time we split a simplex in 2 parts (i.e. maximum depth of
recursive call in \verb#TEST_FACE#).

\begin{center}
\begin{tabular}{|c|c|c|c|c|}
  \hline $\epsilon$ & time & $\mathbb Q$-time& $\mathbb Q$-time/time &max splits \cr
  \hline $5.10^{-1}$ & 0.079s & 0.016s & 20\% &  0 \cr
  \hline $5.10^{-2}$ & 0.150s & 0.022s & 15\% & 6 \cr
  \hline $5.10^{-3}$ & 0.264s & 0.095s & 36\% & 10 \cr
  \hline $5.10^{-4}$ & 0.607s & 0.210s & 35\% & 14 \cr
  \hline $5.10^{-5}$ & 1.638s & 0.656s & 40\% & 17\cr
  \hline $5.10^{-6}$ & 5.380s & 2.548s & 47\% & 20\cr
  \hline $5.10^{-7}$ & 17.087s & 6.732s & 39\% & 24\cr
  \hline
\end{tabular}
\end{center}

The maximum number of splits seems linear in the exponent of $\epsilon$, therefore the
number of subdivision may be at most linear in $\epsilon$ (some simplices needs less
subdivision than others).

We also observe that the final exact test using rational arithmetic never
exceeds half of the total running time. Remark: very small $\epsilon$ would require
multi-precision which we implemented using GMP. But it is far too slow in practice.

\paragraph*{Some curves and surfaces in 2D and 3D}

\begin{figure}
  \begin{center}
    \includegraphics[width=0.4\linewidth,bb=0 0 1194 1035]{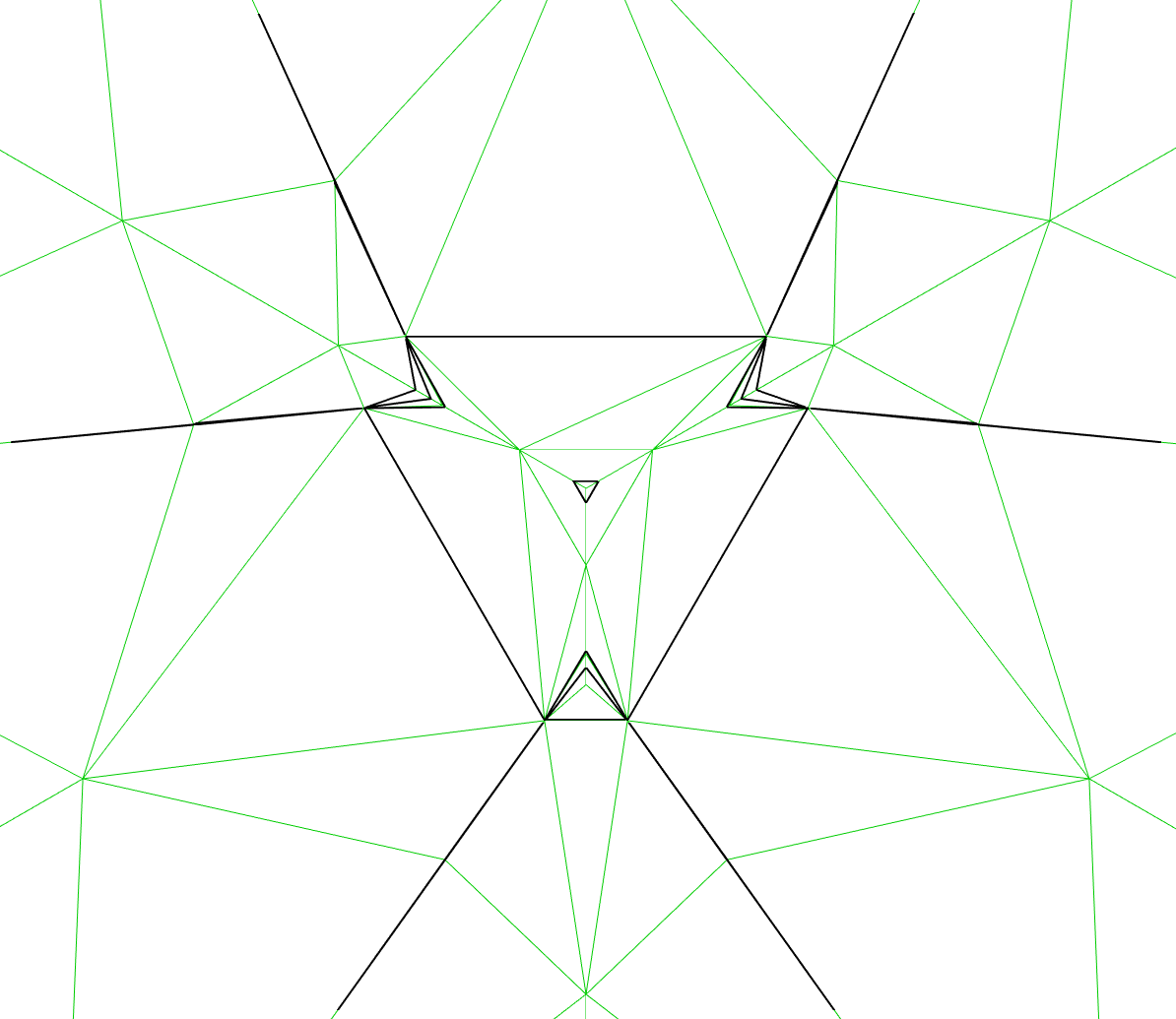}
    \caption{Harnack's sextic}\label{harnack-sextic}
  \end{center}
\end{figure}

\begin{figure}
  \begin{center}
    \includegraphics[height=0.25\paperheight,frame,bb=0 0 360 940]{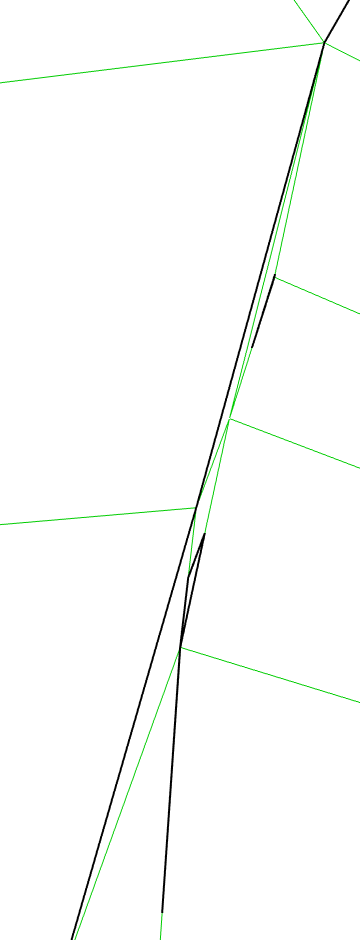}
    \includegraphics[width=0.4\linewidth,bb=0 0 786 1080]{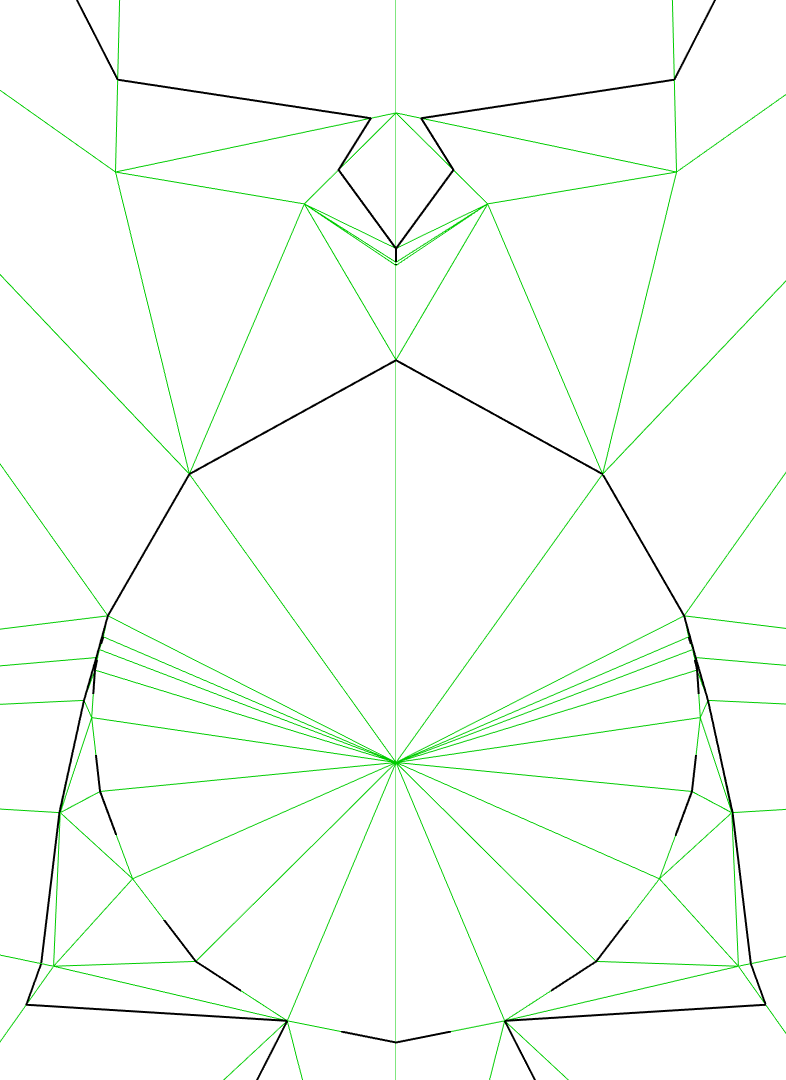}
    \includegraphics[height=0.25\paperheight,frame,bb=0 0 318 556]{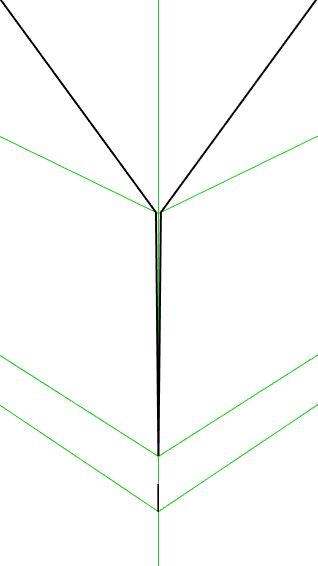}
    \caption{Hilbert's sextix, with two zooms}\label{hilbert-sextic}
  \end{center}
\end{figure}

\begin{figure}
  \begin{center}
    \includegraphics[height=0.25\paperheight,frame,bb=0 0 311 974]{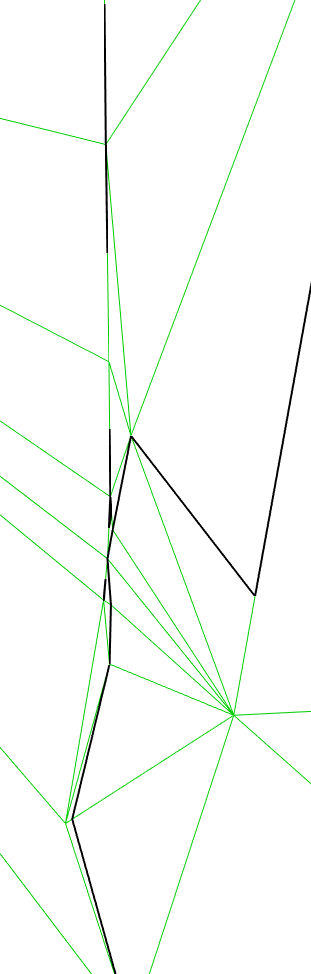}
    \includegraphics[width=0.4\linewidth,bb=0 0 1536 936]{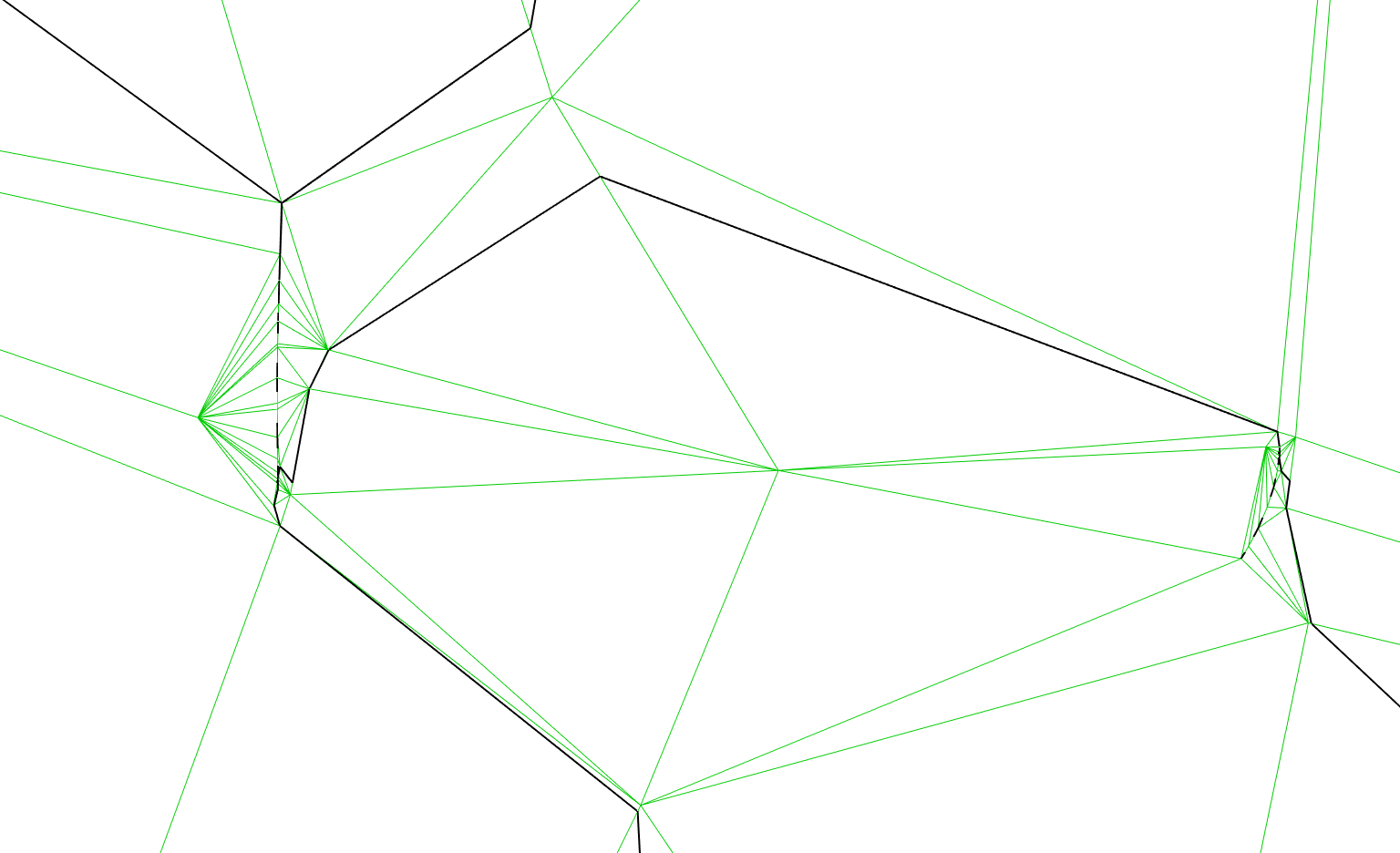}
    \includegraphics[height=0.25\paperheight,frame,bb=0 0 380 1002]{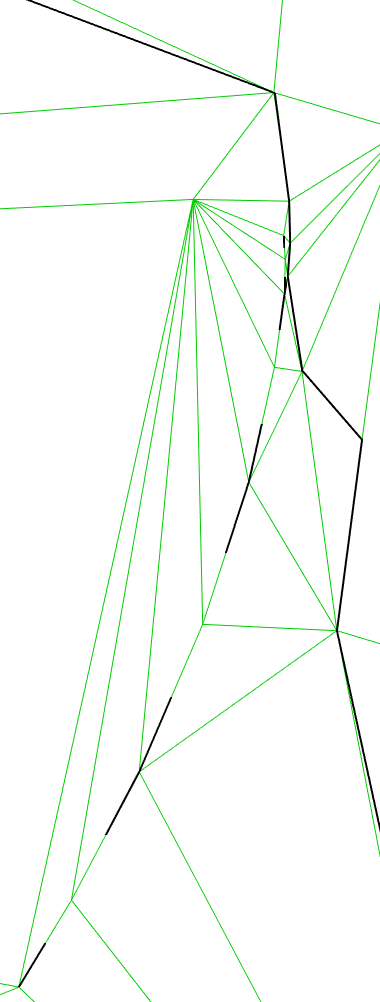}
    \caption{Gudkov's sextix, with two zooms}\label{gudkov-sextic}
  \end{center}
\end{figure}

It is well know that sextic curves have at most 11 components and that this can
be realized in three ways: one oval containing $p$ empty ovals and $10 - p$
empty ovals outside for $p = 1, 5$ or $9$ \cite{Ris93,Vir83}. Construction being
respectively due to Harnak (figure \ref{harnack-sextic}), Hilbert (figure
\ref{hilbert-sextic}) and Gudkov (figure \ref{gudkov-sextic}). Our
implementation succesfully computes the topology of these three curves.

\begin{figure}
  \begin{center}
    \includegraphics[width=0.4\linewidth,bb=0 0 1119 1034]{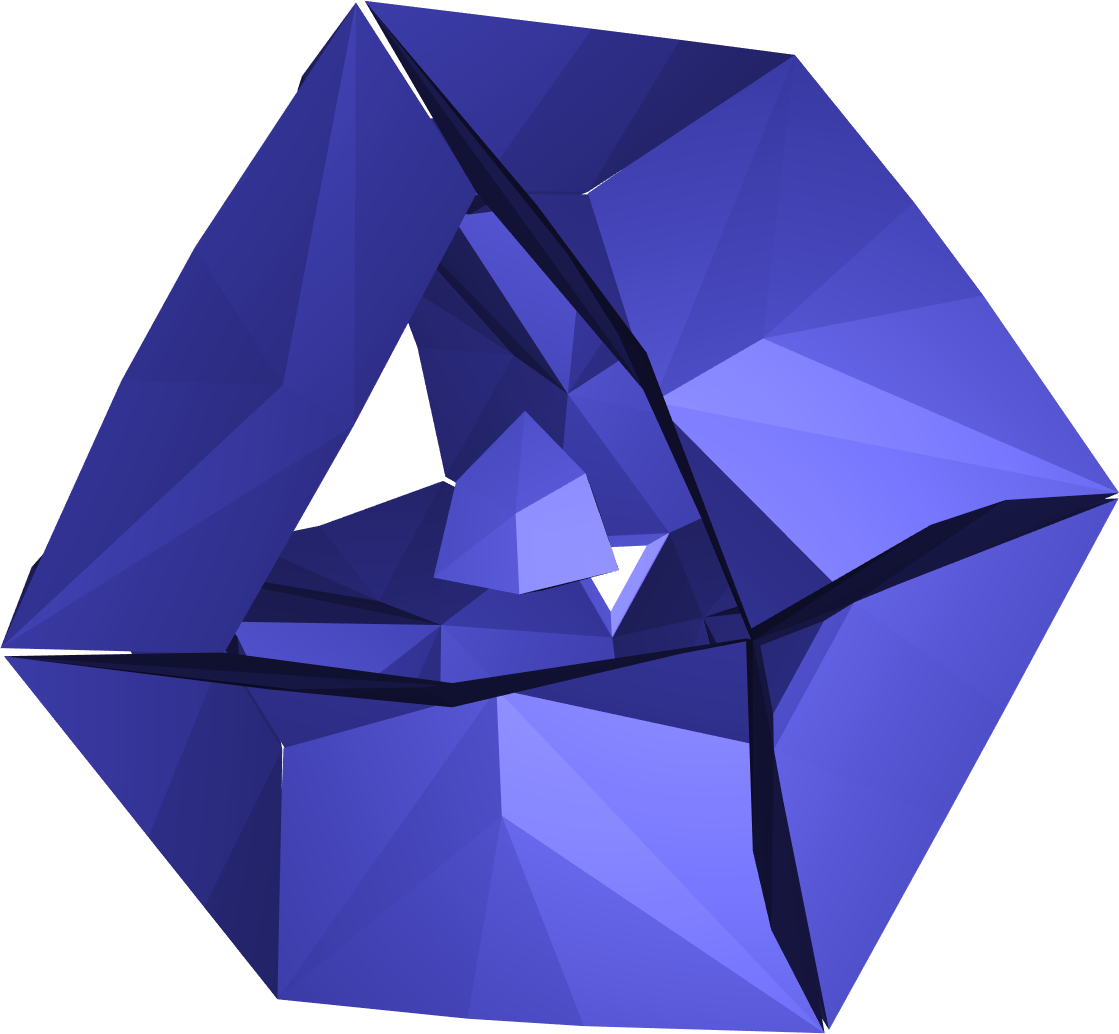}
    \caption{A maximal quartic ``M quartic''}\label{M-quartic}
  \end{center}
\end{figure}

We also experimented succesfullt with two quartic surfaces and two complete
intersection of degree $3\times 2$ and $4 \times 2$:
one maximal quartic (referred as ``M quartic'', figure \ref{M-quartic}) with two components, a sphere and a sphere with 10 handles
and another quartic with 10 spheres (referred as ``M-2 quartic'', figure
\ref{figure:quartic_curve} in the introduction). We also show the intersection of four planes
(product of four linear forms) and a sphere which are used to build the M-2
quartic (referred as ``M-2 quartic $\cap$ S''). Finally, we tested with a 3D curves which is the
intersection of a cone and a cubic surface that gives 5 components (referred as
``cubic $\cap$ cone''). This last example is used in
\cite{Man04,Kol97} to construct Del Pezzo surfaces of degree 1.

Results are summarised in the table below where we give for each example the
codimension (number of polynomials) and projective dimension (number of
variables - 1), the total computing time, the time to check the certificate,
the total number of simplices in the decomposition and the maximum number of
splits  (i.e. maximum depth of
recursive call in \verb#TEST_FACE#).

The number of simplices is smaller than the number of simplices in
the corresponding Newton polytope (counting all quadrants/octants): $4 \times  36 =
144$ for sextic curves and $8 \times 30 = 240$ for quartic surfaces.

\begin{center}
\begin{tabular}{|c|c|c|c|c|c|c|}
  \hline  & codim/dim & time & $\mathbb Q$-time&  simplices & max splits \cr
  \hline Harnack's sextic & 1/2 & 1.430s & 0.142s & 60 &  8 \cr
  \hline Gudkov's sextic  & 1/2 & 15.748s & 1.488s & 76 & 27 \cr
  \hline Hilbert's sextic & 1/2 & 26.036s & 7.724s & 84 & 23 \cr
  \hline M quartic & 1/3 & 5.243s & 0.845s & 204 & 7 \cr
  \hline M-2 quartic & 1/3 & 4.380s & 0.639s & 141 & 7 \cr
  \hline cubic $\cap$ S & 2/3 & 4.050s & 0.280s & 54 & 16 \cr
  \hline M-2 quartic $\cap$ S & 2/3 & 12.809s & 1.304s & 87 & 17 \cr
  \hline
\end{tabular}
\end{center}

\paragraph*{Random varieties}

We also performed experiments with random Polynomials for Bombieri's norm.
This norm is known to give more interesting topology, hence more difficult to
compute than with Euclidien norm. It shows the limit of the current
implementation: we can compute quartic hyper-surfaces up to dimension 5 in
around 20 minutes.

For the zero dimensional case (for which much better approach exists like
msolve), we managed to handle in dimension 4 systems with 3 polynomials of
degree 2 and one of degree 3 (total degree 24) in around 8 minutes (this
takes less than a second with msolve).

You may find in appendix our raw measurements for random polynomials.

It is worth noticing that for all these random tests, the exact test did always
succeed. This is not so surprising as random
varieties are expected to be \emph{smooth enough}.

\paragraph*{Problematic cases}

As mentioned in the previous section, concentric circles (or near to parallel
curves) are currently problematic. We experimented with
$$p(x,y,t) = (x^2+y^2-(1-\alpha)^2)(x^2+y^2-(1+\alpha)^2)$$

We give in the table below the computing time,  number of simplices and max
splits as above. For a difference of radius of $2\alpha=10^{-4}$
we need far more simplices than Newton the polytope for a quartic curve $4 \times 16 =
64$. This is rapidely unfeasible.

\begin{center}
\begin{tabular}{|c|c|c|c|}\hline
  $2\alpha$ & time & simplices & max splits\cr\hline
  1 & 0.083s & 12 & 1 \cr\hline
  $10^{-1}$ & 0.370s & 48 & 6 \cr\hline
  $10^{-2}$ & 1.793s & 128 & 13 \cr\hline
  $10^{-3}$ & 20.476s & 690 & 20\cr\hline
  $10^{-4}$ & 417.259s & 3504 & 28\cr\hline
\end{tabular}
\end{center}


\bibliography{main}
\bibliographystyle{plain}

\appendix
\section{Timings on random polynomials}

Here is how to read a line in these raw results:
\begin{lstlisting}[basicstyle=\small]
  4, (2, 2, 2, 3) => 348.0254s [16.8249s] < 364.8503s (2 samples)
\end{lstlisting}
You can read from left to right:
\begin{itemize}
\item the projective dimension,
\item the degree of each polynomials,
\item the average time to compute a piecewise linear approximation,
\item the standard deviation,
\item the worst time observed and
\item the number of samples.
\end{itemize}

\begin{lstlisting}[basicstyle={\tiny}]
2, (15) => 9.4937s [4.9042s] < 20.0521s (17 samples)
2, (14) => 6.0617s [2.6986s] < 12.4414s (17 samples)
2, (13) => 4.1027s [2.0081s] < 7.2388s (17 samples)
2, (12) => 2.5252s [0.9513s] < 3.8564s (17 samples)
2, (11) => 1.3036s [0.7160s] < 3.0407s (18 samples)
2, (10) => 1.0300s [0.7261s] < 3.2537s (18 samples)
2, (9) => 0.4945s [0.2215s] < 0.8207s (18 samples)
2, (8) => 0.3474s [0.2246s] < 0.9675s (18 samples)
2, (7) => 0.1532s [0.0977s] < 0.4170s (18 samples)
2, (6) => 0.1279s [0.0744s] < 0.2980s (18 samples)
2, (5) => 0.0749s [0.0715s] < 0.2773s (18 samples)
2, (4) => 0.0537s [0.0756s] < 0.3414s (18 samples)
2, (3) => 0.0233s [0.0249s] < 0.0984s (18 samples)
2, (2) => 0.0057s [0.0030s] < 0.0121s (18 samples)
3, (8) => 64.0352s [12.9702s] < 87.0749s (9 samples)
3, (7) => 17.9244s [4.5345s] < 24.4767s (11 samples)
3, (6) => 6.7120s [2.3389s] < 9.9139s (11 samples)
3, (5) => 2.0908s [0.5444s] < 3.0534s (11 samples)
3, (4) => 0.8662s [0.4510s] < 1.6798s (11 samples)
3, (3) => 0.4839s [0.3853s] < 1.4661s (11 samples)
3, (2) => 0.0368s [0.0165s] < 0.0716s (12 samples)
4, (5) => 91.6522s [43.2871s] < 164.7578s (6 samples)
4, (4) => 21.3854s [9.6121s] < 34.8824s (11 samples)
4, (3) => 3.2500s [1.2441s] < 5.7096s (11 samples)
4, (2) => 0.2308s [0.0108s] < 0.2612s (11 samples)
5, (4) => 1119.8743s [183.6737s] < 1303.5480s (2 samples)
5, (3) => 127.2540s [52.1713s] < 192.3819s (3 samples)
5, (2) => 3.8360s [1.4788s] < 5.0987s (3 samples)
6, (3) => 2361.25s [0] < 2361.25s (1 sample)
6, (2) => 14.4625s [0.1260s] < 14.6000s (5 samples)
2, (8, 8) => 5.2380s [3.2701s] < 12.7953s (9 samples)
2, (7, 8) => 6.2327s [2.8787s] < 10.1879s (9 samples)
2, (7, 7) => 2.5434s [1.5762s] < 6.7730s (15 samples)
2, (6, 8) => 2.5903s [1.8971s] < 6.3826s (9 samples)
2, (6, 7) => 2.7967s [1.2523s] < 5.1338s (15 samples)
2, (6, 6) => 2.1022s [2.2025s] < 8.3182s (16 samples)
2, (5, 8) => 1.9484s [0.8972s] < 3.1652s (10 samples)
2, (5, 7) => 1.6489s [1.2579s] < 5.5694s (16 samples)
2, (5, 6) => 2.3762s [2.2781s] < 7.5396s (16 samples)
2, (5, 5) => 1.1877s [0.9082s] < 3.2313s (16 samples)
2, (4, 8) => 1.3558s [0.7226s] < 2.5465s (10 samples)
2, (4, 7) => 1.1043s [0.5974s] < 2.2448s (17 samples)
2, (4, 6) => 0.8182s [0.5225s] < 1.7687s (17 samples)
2, (4, 5) => 0.6955s [0.5440s] < 2.1073s (17 samples)
2, (4, 4) => 0.7081s [0.7159s] < 2.6433s (17 samples)
2, (3, 8) => 1.1122s [0.6825s] < 2.6819s (11 samples)
2, (3, 7) => 1.0376s [0.7008s] < 3.1081s (17 samples)
2, (3, 6) => 0.7728s [0.7388s] < 3.3402s (17 samples)
2, (3, 5) => 0.4749s [0.3775s] < 1.3720s (17 samples)
2, (3, 4) => 0.3867s [0.3490s] < 1.3088s (17 samples)
2, (3, 3) => 0.3086s [0.3875s] < 1.6391s (17 samples)
2, (2, 8) => 0.5736s [0.4633s] < 1.7001s (11 samples)
2, (2, 7) => 0.8804s [1.4395s] < 4.9613s (17 samples)
2, (2, 6) => 0.2350s [0.1518s] < 0.6401s (17 samples)
2, (2, 5) => 0.2294s [0.2023s] < 0.8916s (17 samples)
2, (2, 4) => 0.3407s [0.6465s] < 2.8759s (17 samples)
2, (2, 3) => 0.0939s [0.1056s] < 0.4515s (17 samples)
2, (2, 2) => 0.0578s [0.0621s] < 0.2712s (17 samples)
3, (5, 5) => 88.8346s [39.5227s] < 157.2866s (8 samples)
3, (4, 5) => 43.8346s [13.6107s] < 64.2336s (8 samples)
3, (4, 4) => 23.7975s [14.3489s] < 61.8225s (9 samples)
3, (3, 5) => 28.8840s [25.3741s] < 94.8858s (9 samples)
3, (3, 4) => 13.9328s [6.0994s] < 23.6953s (9 samples)
3, (3, 3) => 6.5742s [2.0125s] < 9.1454s (9 samples)
3, (2, 5) => 16.0180s [13.6073s] < 52.4697s (9 samples)
3, (2, 4) => 11.2471s [12.0677s] < 40.5447s (9 samples)
3, (2, 3) => 1.7853s [1.1207s] < 3.9345s (9 samples)
3, (2, 2) => 0.9715s [1.2720s] < 4.2365s (9 samples)
4, (2, 3) => 120.8822s [43.8194s] < 204.6149s (6 samples)
4, (2, 2) => 52.8344s [44.7051s] < 123.6788s (8 samples)
3, (4, 4, 4) => 57.0216s [24.2671s] < 102.3699s (8 samples)
3, (3, 4, 4) => 63.4938s [22.6517s] < 107.4020s (8 samples)
3, (3, 3, 4) => 39.8249s [18.1909s] < 72.8629s (8 samples)
3, (3, 3, 3) => 21.1542s [13.3550s] < 50.7229s (8 samples)
3, (2, 4, 4) => 52.3019s [55.8282s] < 196.6823s (8 samples)
3, (2, 3, 4) => 14.7853s [9.9172s] < 31.4450s (8 samples)
3, (2, 3, 3) => 11.2110s [4.0091s] < 18.0787s (8 samples)
3, (2, 2, 4) => 14.8802s [9.2282s] < 38.0585s (8 samples)
3, (2, 2, 3) => 5.7705s [2.6911s] < 9.0264s (8 samples)
3, (2, 2, 2) => 5.9955s [4.1095s] < 14.1836s (8 samples)
4, (2, 2, 3) => 480.8930s [111.0565s] < 573.5631s (3 samples)
4, (2, 2, 2) => 166.5942s [86.3898s] < 292.1297s (4 samples)
4, (2, 2, 2, 3) => 348.0254s [16.8249s] < 364.8503s (2 samples)
4, (2, 2, 2, 2) => 247.7984s [24.4447s] < 282.2358s (3 samples)
\end{lstlisting}

\end{document}